\documentclass[11pt,a4paper,reqno]{amsart}
\usepackage[T1]{fontenc}
\usepackage{lmodern}
\usepackage[margin=30mm]{geometry}
\usepackage{amsmath,amssymb,amsthm,mathtools,amscd}
\usepackage{mathrsfs}
\usepackage{tikz-cd}
\usepackage{booktabs,array}
\usepackage{microtype}
\usepackage{xcolor}

\usepackage{enumitem}
\usepackage[colorlinks=true,linkcolor=blue!45!black,citecolor=blue!45!black,urlcolor=blue!45!black]{hyperref}
\usepackage[nameinlink,noabbrev]{cleveref}
\hypersetup{
  pdftitle={A geometric realization of quantum cluster characters and multiplication formula},
  pdfauthor={Chao Shen, Jie Xiao, Fan Xu}
}
\numberwithin{equation}{section}
\newtheorem{theorem}{Theorem}[section]
\newtheorem{introtheorem}{Theorem}

\newtheorem{proposition}[theorem]{Proposition}
\newtheorem{lemma}[theorem]{Lemma}
\newtheorem{corollary}[theorem]{Corollary}
\theoremstyle{definition}
\newtheorem{definition}[theorem]{Definition}
\newtheorem{remark}[theorem]{Remark}
\newtheorem{example}[theorem]{Example}
\newcommand{\Q}{{\overline{\mathbb Q}_{\ell}}}
\newcommand{\Fq}{\mathbb F_q}
\newcommand{\Fr}{\operatorname{Fr}_q}
\newcommand{\M}{\mathfrak M}
\newcommand{\Db}{D^b}

\newcommand{\Cq}{\mathscr C_q}
\newcommand{\HH}{H}
\newcommand{\TT}{\mathbb T}
\newcommand{\KK}{K_0}
\newcommand{\eps}{\varepsilon}
\newcommand{\ii}{\mathbf i}
\newcommand{\iiO}{\mathbf i_0}
\newcommand{\Y}{\mathcal Y}

\newcommand{\bal}[1]{[#1]\left(\frac{#1}{2}\right)}
\newcommand{\RGC}{\mathrm R\Gamma_c}

\newcommand{\boxt}{\boxtimes}
\newcommand{\eu}[2]{\langle #1,#2\rangle}
\newcommand{\sy}[2]{(#1,#2)}
\newcommand{\ef}[2]{\mathsf E_{\ii}(#1,#2)}

\newcommand{\wt}[1]{|#1|_{\ii}}
\newcommand{\ellone}[1]{|#1|_1}
\DeclareMathOperator{\Ext}{Ext}
\DeclareMathOperator{\Hom}{Hom}
\DeclareMathOperator{\Aut}{Aut}

\DeclareMathOperator{\Gr}{Gr}
\DeclareMathOperator{\Filt}{Filt}
\DeclareMathOperator{\Rep}{Rep}
\DeclareMathOperator{\GL}{GL}
\DeclareMathOperator{\Tr}{Tr}

\DeclareMathOperator{\pt}{pt}

\DeclareMathOperator{\Spec}{Spec}

\setlist[itemize]{leftmargin=*,itemsep=2pt,topsep=4pt}
\setlist[enumerate]{leftmargin=*,itemsep=2pt,topsep=4pt}
\title[Quantum cluster characters and multiplication formula]{A geometric realization of quantum cluster characters and multiplication formula}
\author{Chao Shen}
\address{School of Mathematical Sciences,\newline
  \hspace*{\parindent}University of Science and Technology of China, Hefei 230026, China}
\email{sc26@ustc.edu.cn}
\author{Jie Xiao}
\address{School of Mathematical Sciences,\newline
  \hspace*{\parindent}Beijing Normal University, Beijing 100875, China}
\email{jxiao@bnu.edu.cn}
\author{Fan Xu}
\thanks{Chao Shen is partially supported by the National Natural Science Foundation of China (NSFC, Grant No. 12471068). Jie Xiao is partially supported by NSFC (Grant No. 12471030). Fan Xu is partially supported by NSFC (Grant No. 12371036).}
\address{Department of Mathematical Sciences,\newline
  \hspace*{\parindent}Tsinghua University, Beijing 100084, China}
\email{fanxu@tsinghua.edu.cn}
\date{}

\begin{document}

\begin{abstract}
We construct a character functor for finite acyclic quivers over finite
fields by pairing mixed complexes on representation stacks with Lusztig
sheaves. Lusztig's restriction formula gives a multiplicative filtration
for this functor. Taking Frobenius traces recovers the Hall
characters of Berenstein and Rupel and, as a consequence,
Rupel's quantum cluster characters.
We also construct distinguished triangles for extension families and use
finite filtrations to relate their middle terms to induction products,
which correspond to Hall convolution under Frobenius trace.
Applying the character functor to these triangles gives a categorical
realization of the general multiplication formula.
\end{abstract}

\maketitle
\enlargethispage{3pt}
\tableofcontents

\section{Introduction}
\label{sec:introduction}

\subsection{Background and motivation}

Cluster algebras are commutative algebras generated by cluster variables,
which are grouped into overlapping clusters related by mutation
\cite{FominZelevinsky2002}. Since mutation defines these variables
recursively, a basic problem is to find their explicit Laurent expansions
with respect to a fixed initial seed. The search for such formulas is
closely tied to additive categorification via cluster categories and
their generalizations \cite{BMRRT2006,Amiot2009,KellerYang2011}.
For Dynkin quivers, Caldero and
Chapoton obtained such expressions
using Euler characteristics of quiver Grassmannians \cite{CC2006}.
Caldero and Keller extended this construction to acyclic quivers
\cite{CalderoKeller2006}. Derksen, Weyman and Zelevinsky gave analogous
cluster character formulas using quivers with potentials
\cite{DWZ2008,DWZ2010}.
Palu and Plamondon further developed cluster characters for
$2$-Calabi--Yau triangulated categories
\cite{Palu2008,Plamondon2011}.
More general multiplication formulas express products of characters
through extension data \cite{XiaoXu2009,XiaoXu2010,Xu2010,Palu2012}.
Quiver Grassmannians enter these cluster characters through
their Euler characteristics.

The Hall-algebra viewpoint is already present in the work of Caldero
and Chapoton. In finite type, Hall polynomials relate finite-field point
counts to the Euler characteristics occurring in their formula,
suggesting a counting refinement \cite{CC2006}.
Qin \cite{Qin2012} and Rupel \cite{Rupel2015} use finite-field methods
to establish quantum cluster character formulas for acyclic quivers
and acyclic valued quivers, respectively. They establish counting
polynomials for quiver Grassmannians of rigid representations,
yielding formulas in a formal quantum parameter.
Berenstein and Rupel use flag counts to construct algebra homomorphisms
from dual Hall algebras, recovering quantum cluster characters
for doubled repetition-free source-adapted words
\cite{BR2015}.
The Hall-algebra approach adds an associative product defined by
extension counts to the additive categorical picture, and the character
maps respect this product \cite{BR2015,ChenDingXu2016,DingXuZhang2020,FuPengZhang2023}.
In these finite-field formulas, the geometry of flags and
extensions enters through point counts.

Motivic Hall algebras retain more geometric information by working with
classes of stacks rather than point counts \cite{Bridgeland2012}. Nagao uses these algebras
to study cluster transformations \cite{Nagao2013}, and Bridgeland
uses them to construct scattering diagrams
\cite{Bridgeland2017}. Davison lifts wall-crossing identities to
isomorphisms in categories of monodromic mixed Hodge modules
\cite{Davison2018}.

Of course, the above list of references on classical and quantum
cluster characters is far from complete.
Inspired by Davison's sheaf-theoretic approach, we seek to lift the finite-field Hall character
homomorphisms and multiplication formulas to mixed complexes.

\subsection{Main results}

We construct a character functor on mixed complexes over representation
stacks of finite acyclic quivers. Our main results describe its
multiplicative filtration, an extension triangle, and the coefficients
of the boundary term. Taking Frobenius traces recovers the Hall
character homomorphisms and the general multiplication formula.

Let $Q$ be a finite acyclic quiver over $k=\Fq$, and let $\M_\alpha$
be its representation stack in dimension $\alpha$. A word
$\ii=(i_1,\ldots,i_m)$ and a tuple $\mathbf a\in\mathbb N^m$
specify flags with successive quotients $S_{i_k}^{\oplus a_k}$.
The proper flag projection defines
$\mathcal K_{\mathbf a}=\pi_{\mathbf a!}\Q$. We use the
compact-support pairing
\[
 \mathsf P_\alpha(\mathcal F,\mathcal K)
   =\RGC(\M_\alpha,\mathcal F\otimes\mathcal K).
\]
For the point object $A_L=x_{L!}\Q$ associated with a representation
$L$, this gives
\[
 \mathsf P_\alpha(A_L,\mathcal K_{\mathbf a})
   \simeq\RGC(\Filt_{\ii,\mathbf a}(L),\Q).
\]
These pairings, with the normalization of
Definition~\ref{geo:character-definition}, define the character functor
$\Y_\ii$. Write $\star=p_!q^!$ for induction and $\odot$
for the twisted product of graded coefficient complexes in
\eqref{geo:target-product}.

\begin{introtheorem}[Multiplicative filtration, see Theorem~\ref{geo:main-theorem}]
\label{intro:theorem-a}
Let $\mathcal F\in\Db_m(\M_\alpha)$ and
$\mathcal G\in\Db_m(\M_\beta)$ be bounded mixed complexes.
For each flag type $\mathbf c$ of weight $\alpha+\beta$, the complex
$\Y_\ii(\mathcal F\star\mathcal G)_{\mathbf c}$ admits a finite
filtration, functorial in $\mathcal F$ and $\mathcal G$, such that
\[
 \operatorname{gr}\!\left(
   \Y_\ii(\mathcal F\star\mathcal G)_{\mathbf c}\right)
   \simeq
   \bigl(\Y_\ii(\mathcal F)\odot\Y_\ii(\mathcal G)\bigr)_{\mathbf c}.
\]
Here $\operatorname{gr}$ denotes the direct sum of the filtration layers.
\end{introtheorem}

The proof uses
Lusztig's restriction filtration and the pairing above. In Grothendieck
groups, the class of a filtered object is the sum of the classes of its
layers, so the character functor induces a ring homomorphism.

Taking Frobenius traces then gives an algebra
homomorphism to the quantum torus $\TT_\ii$ associated with the word
$\ii$. Its based-monomial coefficients agree with those of
Berenstein and Rupel \cite{BR2015}, as shown in
Corollary~\ref{char:BR}. Remark~\ref{char:comparison} compares the products.
This construction also recovers Reineke's integration map \cite{Reineke2006}
and the Feigin homomorphism \cite{Reineke2001}.
Let $\iiO$ be a repetition-free source-adapted word on all vertices.
For the doubled word $\ii=(\iiO,\iiO)$, a change
of coordinates gives the quantum cluster character
(Corollary~\ref{cc:character}). The point objects of indecomposable rigid representations thus realize all non-initial quantum cluster variables (Corollary~\ref{cc:rigid-specialization}).

For representations $M,N$ of dimension vectors $m,n$, put
\[
 V=\Ext^1_Q(M,N),\qquad r=\dim V>0.
\]
The quotient stack $[V/\mathbb G_m]$ carries
a family of middle terms $\phi:[V/\mathbb G_m]\to\M_{m+n}$ whose
restriction to the open immersion
$j:\mathbb P(V)\hookrightarrow[V/\mathbb G_m]$ parametrizes
nonsplit extensions.

\begin{introtheorem}[Extension triangle, see Theorem~\ref{cc:sheaf-exchange}]
\label{intro:theorem-b}
For arbitrary $M,N$ with $r>0$, there is a complex
$\mathcal H_{M,N}$ and a distinguished triangle
\begin{equation}\label{eq:intro-general-triangle}
 \mathcal E_{M,N}
 \longrightarrow\mathcal H_{M,N}
 \longrightarrow\mathcal B_{M,N}[2r](r)
 \longrightarrow,
\end{equation}
where
\[
 \mathcal E_{M,N}=(\phi j)_!\Q,\qquad
 \mathcal B_{M,N}=\phi_!Rj_*\Q[-1].
\]
The middle term has a finite filtration with layers
\[
 (A_M\star A_N)\otimes\Q[-2s](-s),\qquad 0\leq s<r.
\]
\end{introtheorem}

The next result determines the
Grassmannian coefficients of the boundary term $\mathcal B_{M,N}$.

Take $\ii=(\iiO,\iiO)$ and put $\alpha=m+n$.
The Lusztig sheaf $\mathcal K_{(\alpha-e,e)}$ then corresponds to
the relative quiver Grassmannian of subrepresentations of dimension $e$.
For subrepresentations
$U\subset M$ and $W\subset N$, pullback along $U\hookrightarrow M$
and pushout along $N\twoheadrightarrow N/W$ give the map
\[
 \beta_{U,W}:V\longrightarrow\Ext^1_Q(U,N/W).
\]
It is surjective since $\Rep_k(Q)$ is hereditary.
Over $S_{u,v}=\Gr_u(M)\times\Gr_v(N)$, define
\[
 \Xi_{u,v}=\{(U,W,[\eta]):
   [\eta]\in\mathbb P(\operatorname{im}D_k\beta_{U,W})\}
   \subset S_{u,v}\times\mathbb P(D_kV),
\]
where $D_k=\Hom_k(-,k)$. The closed scheme structure is described in
Subsection~\ref{bd:sec:lifting}.

\begin{samepage}
\begin{introtheorem}[Boundary coefficients, see Theorem~\ref{bd:thm:main}]
\label{intro:theorem-c}
For every dimension vector $0\leq e\leq\alpha$, the complex
$\mathsf P(\mathcal B_{M,N},\mathcal K_{(\alpha-e,e)})$ admits a finite
filtration with layers
\[
 \RGC(\Xi_{u,v},\Q)[-2t_{u,v}](-t_{u,v}),
 \qquad t_{u,v}=r+\eu{u}{n-v},\quad u+v=e,
\]
indexed by $u\leq m$, $v\leq n$, and $u+v=e$.
Here $\langle-,-\rangle$ is the Euler form of $Q$.
\end{introtheorem}
\end{samepage}

Applying the character functor to the triangle in Theorem~\ref{intro:theorem-b}
and taking Frobenius traces gives a formula involving the boundary character.
Theorem~\ref{intro:theorem-c} and Proposition~\ref{gm:ordered} identify
this term and yield the multiplication formula for ordinary acyclic quivers with principal
coefficients (Theorem~\ref{gm:general-multiplication}).
The factor
$[r]_q=1+q+\cdots+q^{r-1}$ comes from the middle filtration in
Theorem~\ref{intro:theorem-b}. For rigid exchange pairs,
Corollary~\ref{gm:exchange-mutation} recovers the quantum exchange
relation for the characters of Corollary~\ref{cc:rigid-specialization}.
For algebraic approaches to cluster multiplication formulas, see
\cite{ChenDingZhang2023,ChenXiaoXu2023}.

Our contribution is the realization of the formula in
Theorem~\ref{gm:general-multiplication}
by an extension triangle and
boundary coefficient complexes, together with the multiplicative
filtration connecting these constructions to quantum characters.

\smallskip
\noindent\textit{Conventions and references.}\enspace
We originally formulated our constructions in Lusztig's equivariant
setting \cite{Lusztig1991,Lusztig1993}.
We use quotient stacks here so that the relation with Hall
convolution is explicit. In this notation, Lusztig's induction correspondence
is the correspondence of short exact sequences in \eqref{eq:hall-correspondence}.
For background on sheaves on algebraic stacks, see \cite{LO2008}.
For quiver representations and Hall algebras, see
\cite{AssemSimsonSkowronski2006} and \cite{Ringel1996}, respectively.

\subsection{Organization of the paper}

In Section~\ref{sec:conventions}, we introduce the representation stacks and
pairing. In Section~\ref{sec:flags}, we construct the character functor
and use Lusztig's restriction filtration to prove Theorem~\ref{intro:theorem-a}.
In Section~\ref{sec:characters}, we identify
the resulting character homomorphisms.
In Section~\ref{cc:sheaf-triangle}, we construct the extension triangle.
In Section~\ref{sec:boundary-multiplication}, we determine its boundary
coefficients and derive the multiplication formula.
\section{Sheaf-theoretic framework on representation stacks}
\label{sec:conventions}

\subsection{Quivers and representation stacks}
Let $Q=(Q_0,Q_1)$ be a finite quiver without oriented cycles, and let
$k=\Fq$ and $\pt=\Spec k$. For $h\in Q_1$, write $s(h)$ and $t(h)$ for its source and target.
For dimension vectors $\alpha,\beta\in\mathbb N^{Q_0}$, set
\begin{equation}\label{eq:euler}
 \eu{\alpha}{\beta}
 =\sum_{i\in Q_0}\alpha_i\beta_i
       -\sum_{h\in Q_1}\alpha_{s(h)}\beta_{t(h)},
 \qquad
 \sy{\alpha}{\beta}=\eu{\alpha}{\beta}+\eu{\beta}{\alpha}.
\end{equation}

The category $\Rep_k(Q)$ is hereditary. For $M,N\in\Rep_k(Q)$,
the Euler form satisfies
\[
 \eu{\dim M}{\dim N}
 =\dim_k\Hom_Q(M,N)-\dim_k\Ext^1_Q(M,N).
\]

Let $S_i$ be the simple representation at vertex $i$,
and let $\eps_i$ be its dimension vector. Put
$\ellone{\alpha}=\sum_i\alpha_i$.

Choose graded vector spaces of dimension $\alpha$ and put
\[
 E_\alpha=\bigoplus_{h\in Q_1}
       \Hom(k^{\alpha_{s(h)}},k^{\alpha_{t(h)}}),\qquad
 G_\alpha=\prod_{i\in Q_0}\GL_{\alpha_i},\qquad
 \M_\alpha=[E_\alpha/G_\alpha].
\]
The stack $\M_\alpha$ is smooth of dimension $-\eu{\alpha}{\alpha}$.

Let $\mathfrak E_{\alpha,\beta}$ be the algebraic stack of short exact
sequences of $Q$-representations
\begin{equation}\label{eq:extension-convention}
 0\longrightarrow N\longrightarrow L\longrightarrow M\longrightarrow0,
 \qquad \dim M=\alpha,\quad\dim N=\beta.
\end{equation}
The Hall correspondence is
\begin{equation}\label{eq:hall-correspondence}
 \M_\alpha\times\M_\beta
 \xleftarrow{\ q\ }\mathfrak E_{\alpha,\beta}
 \xrightarrow{\ p\ }\M_{\alpha+\beta}.
\end{equation}

The maps are defined by
\[
 \begin{aligned}
 q(0\to N\to L\to M\to0)&=(M,N),\\
 p(0\to N\to L\to M\to0)&=L.
 \end{aligned}
\]
The morphism $p$ is representable and proper, with quiver Grassmannians
as geometric fibres. The morphism $q$ is a vector bundle stack of rank
$-\eu{\alpha}{\beta}$. The diagram \eqref{eq:hall-correspondence}
is Lusztig's induction correspondence written in quotient-stack
notation \cite{Lusztig1991,Lusztig1993}.

\subsection{Mixed complexes, twists, and convergence}

Fix a prime $\ell\ne\operatorname{char}k$ and an embedding
$\iota:\overline{\mathbb Q}_\ell\hookrightarrow\mathbb C$.
For an algebraic stack $\mathfrak X$ of finite type over $k$, write
$\Db_m(\mathfrak X)$ for the bounded derived category of constructible
mixed $\Q$-complexes on $\mathfrak X$.
For background on mixed complexes and weights over finite fields, see
\cite{BBD1982}.
Fix an algebraic closure
$\bar k$ of $k=\Fq$. For a complex ${\mathcal F}$ on $\mathfrak X$, write
${\mathcal F}_{\bar k}$ for its pullback to $\mathfrak X_{\bar k}$, equipped
with the induced Frobenius structure. Write $\Db_c(\mathfrak X_{\bar k})$
for the bounded constructible derived category on the geometric stack
$\mathfrak X_{\bar k}$. We use the six operations for
algebraic stacks \cite{LO2008,Sun2012}. Tensor products and all functors
are derived.

For a bounded mixed complex $C$ on $\pt$, its Frobenius trace is
\[
 \Tr(\Fr,C)
 =\sum_j(-1)^j\iota\!\left(\Tr(\Fr;H^j(C_{\bar k}))\right),
\]
where $\Fr$ denotes geometric Frobenius. For
${\mathcal F}\in\Db_m(\M_\alpha)$, the associated trace function is
$f_{\mathcal F}(L)=\Tr(\Fr,x_L^*{\mathcal F})$, where $x_L:\pt\to\M_\alpha$ is the point
defined by $L\in\Rep_k(Q)$ of dimension $\alpha$.

The symbol $(1)$ denotes the Tate object on which geometric Frobenius acts
by $q^{-1}$. Choose its square root $(1/2)$ so that its Frobenius eigenvalue
has image $q^{-1/2}>0$ under $\iota$. We use paired shifts
$C[d](d/2)$, with $d\in\mathbb Z$, and set
\begin{equation}\label{eq:paired-shift}
 w=-q^{1/2}.
\end{equation}
The shift and Tate twist act on Frobenius traces by
\begin{equation}\label{eq:trace-shift}
 \Tr(\Fr,C[d](r))=(-1)^d q^{-r}\Tr(\Fr,C),\qquad
 \Tr(\Fr,C\bal{d})=w^{-d}\Tr(\Fr,C).
\end{equation}
We use $w$ as the numerical Hall and quantum torus parameter.
The minus sign records the cohomological shift in $\bal{1}$.

\begin{example}[$B\mathbb G_m$]\label{ex:BGm}
The cohomology ring of $B\mathbb G_m$ is $\Q[c_1]$, with
$H^{2m}(B\mathbb G_{m,\bar k},\Q)=\Q(-m)$
for $m\ge0$ and no odd cohomology. Since
$\dim B\mathbb G_m=-1$, Poincar\'e duality gives
\[
 H_c^{-2-2m}(B\mathbb G_{m,\bar k},\Q)
 \simeq H^{2m}(B\mathbb G_{m,\bar k},\Q)^\vee(1)
 \simeq\Q(m+1).
\]
Thus $\Q(n)$ occurs in degree $-2n$ for each $n\ge1$.
Its Frobenius eigenvalue is $q^{-n}$, and the degree is even, so
\[
 \Tr(\Fr,\RGC(B\mathbb G_m,\Q))
   =\sum_{n\ge1}q^{-n}=\frac1{q-1}.
\]
See \cite[Example~7.1]{Sun2012}. This agrees with the weighted
point count, since every one-dimensional $k$-vector space is
isomorphic to $k$:
\[
 \#B\mathbb G_m(k)
 =\sum_{\substack{[L]\\\dim_k L=1}}\frac{1}{\#\Aut_k(L)}
 =\frac{1}{\#\GL_1(k)}=\frac{1}{q-1}.
\]
\end{example}

Thus compactly supported cohomology on a stack need not be bounded,
even when its Frobenius trace converges. The coefficient complexes
of our character functors lie in the
following category of convergent mixed complexes on $\pt$.

\begin{definition}\label{def:convergent}
Let $\Cq$ be the full subcategory of $D^-_m(\pt)$, the bounded-above
mixed derived category on $\pt$, consisting of complexes $C$ with
finite-dimensional geometric cohomology in each degree such that
\begin{equation}\label{eq:convergence}
 \sum_{j\in\mathbb Z}\ \sum_{\lambda\in\operatorname{Eig}(\Fr;H^j(C_{\bar k}))}
       |\iota(\lambda)|^s<\infty\qquad\text{for every }s\in\mathbb R_{>0}.
\end{equation}

Eigenvalues are counted with algebraic multiplicity. The Frobenius trace
of $C$ is defined by the absolutely convergent sum
\[
 \Tr(\Fr,C)=\sum_{j\in\mathbb Z}(-1)^j
   \iota\!\left(\Tr(\Fr;H^j(C_{\bar k}))\right).
\]

\end{definition}

This is the convergence condition of \cite[Definition~4.1]{Sun2012}.
In that formalism we use stratifiable complexes whenever a direct image is
taken. All inputs on representation and flag stacks in this paper are
bounded constructible mixed complexes, hence are stratifiable
by \cite[Lemma~3.4(v)]{Sun2012}.

\begin{lemma}\label{lem:convergent}

The category $\Cq$ has the following properties.
\begin{enumerate}[label=\textup{(\roman*)}]
\item For $U,V\in\Cq$, $n\in\mathbb Z$, and $r\in\frac12\mathbb Z$,
\[
 U\oplus V,\qquad U[n](r),\qquad U\otimes V\quad\text{belong to }\Cq.
\]
Their Frobenius traces satisfy
\[
 \begin{aligned}
 \Tr(\Fr,U\oplus V)&=\Tr(\Fr,U)+\Tr(\Fr,V),\\
 \Tr(\Fr,U[n](r))&=(-1)^nq^{-r}\Tr(\Fr,U),\\
 \Tr(\Fr,U\otimes V)&=\Tr(\Fr,U)\Tr(\Fr,V).
 \end{aligned}
\]
\item For any morphism $u:U\to V$ in $\Cq$, the mapping cone
$W=\operatorname{Cone}(u)$ also belongs to $\Cq$. It fits into a
distinguished triangle
\[
 U\xrightarrow{\ u\ }V\longrightarrow W\longrightarrow U[1],
 \qquad
 \Tr(\Fr,V)=\Tr(\Fr,U)+\Tr(\Fr,W).
\]
Thus $\Cq$ is a full triangulated subcategory of $D^-_m(\pt)$, and
Frobenius trace is additive on its distinguished triangles.
\item Let $\mathfrak X$ be an algebraic stack of finite type over $k$.
For every $\mathcal F\in\Db_m(\mathfrak X)$,
\[
 \RGC(\mathfrak X,\mathcal F)\in\Cq.
\]
\end{enumerate}

\end{lemma}
\begin{proof}
The closure under cones and additivity of traces are
\cite[Lemma~4.3]{Sun2012}. For tensor products, write $S_s(C)$ for
the sum in \eqref{eq:convergence}. K\"unneth gives
$S_s(U\otimes V)=S_s(U)S_s(V)$. Each cohomological degree involves
only finitely many summands because both complexes are bounded above.
Absolute convergence gives the trace identity. Direct sums, shifts and
twists have the stated effect on the eigenvalues. Finally, bounded
constructible mixed inputs are stratifiable and convergent, so
\cite[Theorem~4.2(i)]{Sun2012}, applied to the structure morphism,
proves \textup{(iii)}.
\end{proof}

Throughout the paper, $\KK(\mathscr D)$ denotes the Grothendieck group of
a triangulated category $\mathscr D$, generated by isomorphism classes
$[A]$ of objects, with relations $[B]=[A]+[C]$ for distinguished triangles
$A\to B\to C\to A[1]$.
For a finite filtration $0=E_0\to E_1\to\cdots\to E_s=E$ in
$\Cq$, the layers $L_j$ are the cones in the distinguished triangles
\[
 E_{j-1}\longrightarrow E_j\longrightarrow L_j\longrightarrow E_{j-1}[1].
\]
They play the role of successive quotients. The triangle relations give
\[
 [E]=\sum_{j=1}^s[L_j]\quad\text{in }\KK(\Cq),\qquad
 \Tr(\Fr,E)=\sum_{j=1}^s\Tr(\Fr,L_j).
\]
In particular, $[\mathcal F[d]]=(-1)^d[\mathcal F]$.

For $\mathcal F\in\Db_m(\M_\alpha)$, the stack trace formula
\cite[Theorem~4.2(ii)]{Sun2012} gives
\begin{equation}\label{eq:stack-trace}
 \begin{aligned}
 \Tr(\Fr,\RGC(\M_\alpha,\mathcal F))
 &=\sum_{\substack{[L]\\\dim L=\alpha}}
   \frac{f_{\mathcal F}(L)}{a_L}\\
 &=\frac{1}{|G_\alpha(k)|}\sum_{x\in E_\alpha(k)}f_{\mathcal F}(L_x).
 \end{aligned}
\end{equation}
Here $[L]$ runs over the $k$-isomorphism classes of representations of
dimension $\alpha$, $a_L=|\Aut(L)(k)|$, and $L_x$ is the representation
defined by $x\in E_\alpha(k)$. The second equality is the
orbit--stabilizer formula.

\subsection{Induction and restriction}
Let $\mathscr D_Q=\bigoplus_\alpha\Db_m(\M_\alpha)$, with finite support in
the dimension vector. For $\mathcal F\in\Db_m(\M_\alpha)$ and
$\mathcal G\in\Db_m(\M_\beta)$, define induction using the Hall
correspondence \eqref{eq:hall-correspondence} by
\begin{equation}\label{eq:raw-convolution}
 \mathcal F\diamond\mathcal G=p_!q^*(\mathcal F\boxt\mathcal G).
\end{equation}
Writing $\mathbb D$ for Verdier duality, define the dual product by
\begin{equation}\label{eq:star}
 \mathcal F\star\mathcal G
 =\mathbb D\bigl((\mathbb D\mathcal F)\diamond(\mathbb D\mathcal G)\bigr)
 \simeq p_!q^!(\mathcal F\boxt\mathcal G)
 \simeq(\mathcal F\diamond\mathcal G)[-2\eu{\alpha}{\beta}]
                                     (-\eu{\alpha}{\beta}).
\end{equation}
Here the isomorphisms follow from properness of $p$ and smooth purity
for $q$, whose relative dimension is $-\eu{\alpha}{\beta}$.
Both functors preserve bounded mixed complexes.

Base change for the stack of flags $0\subset L_1\subset L_2\subset L_3=L$
gives
\[
 a_{\mathcal F,\mathcal G,\mathcal H}:
 (\mathcal F\diamond\mathcal G)\diamond\mathcal H
 \xrightarrow{\sim}
 \mathcal F\diamond(\mathcal G\diamond\mathcal H).
\]
For $\star$, bilinearity of the Euler form gives
\[
 \eu{\alpha}{\beta}+\eu{\alpha+\beta}{\gamma}
 =\eu{\beta}{\gamma}+\eu{\alpha}{\beta+\gamma}.
\]
Thus the two iterated products have the same total paired shift.
We fix associativity constraints compatible with the unit, making
$(\mathscr D_Q,\star)$ a monoidal category. In particular, $\star$
induces an associative product on $K_0(\mathscr D_Q)$.
With $\mathbf1=\Q_{\M_0}$, the unit isomorphisms are
\[
 \mathbf1\diamond\mathcal F\simeq\mathcal F\diamond\mathbf1\simeq\mathcal F,
 \qquad
 \mathbf1\star\mathcal F\simeq\mathcal F\star\mathbf1\simeq\mathcal F.
\]
Using the same correspondence \eqref{eq:hall-correspondence} in the
opposite direction, define restriction by
\[
 \mathcal R_{\alpha,\beta}(\mathcal K)=q_!p^*\mathcal K,
 \qquad \mathcal K\in\Db_m(\M_{\alpha+\beta}).
\]
Here $q_!$ preserves bounded mixed complexes, since $q$ is a
vector-bundle stack.

\subsection{A Lusztig-type pairing and adjunction}
\label{geo:pairing-section}

For $\mathcal F,\mathcal K\in\Db_m(\M_\alpha)$ define
\begin{equation}\label{geo:pairing-definition}
 \mathsf P_\alpha(\mathcal F,\mathcal K)
       =\RGC(\M_\alpha,\mathcal F\otimes \mathcal K).
\end{equation}
The pairing takes values in $\Cq$ by Lemma~\ref{lem:convergent}.
It preserves distinguished triangles in each variable, since derived
tensor product and $\RGC$ do so. On a product of representation stacks
we use the notation $\mathsf P_{\alpha,\beta}$ for the same
compact-support pairing. This is a complex-valued quotient-stack variant of Lusztig's
equivariant pairing \cite{Lusztig1993}.

The following formula is due to Lusztig
\cite[Lemma~9.2.9 and \S9.2.10]{Lusztig1993}.

\begin{lemma}[Adjunction]
\label{geo:pairing-transposition}
Let $\mathcal F\in\Db_m(\M_\alpha)$ and
$\mathcal G\in\Db_m(\M_\beta)$. For every
$\mathcal K\in\Db_m(\M_{\alpha+\beta})$ there is a natural isomorphism
\begin{equation}\label{geo:pairing-transposition-equation}
 \mathsf P_{\alpha+\beta}(\mathcal F\diamond\mathcal G,\mathcal K)
 \simeq
 \mathsf P_{\alpha,\beta}
  (\mathcal F\boxt\mathcal G,\mathcal R_{\alpha,\beta}\mathcal K).
\end{equation}
For a representation point $x_L:\pt\to\M_\alpha$
and $A_L=x_{L!}\Q$, one has, for every $\mathcal K\in\Db_m(\M_\alpha)$,
\begin{equation}\label{geo:point-evaluation}
 \mathsf P_\alpha(A_L,\mathcal K)\simeq x_L^*\mathcal K.
\end{equation}
\end{lemma}

\begin{proof}
Set $\mathcal H=\mathcal F\boxt\mathcal G$.
The projection formula and
composition of compact-support pushforwards give
\begin{align*}
 \RGC(\M_{\alpha+\beta},p_!q^*\mathcal H\otimes \mathcal K)
 &\simeq\RGC(\mathfrak E_{\alpha,\beta},q^*\mathcal H\otimes p^*\mathcal K)\\
 &\simeq\RGC(\M_\alpha\times\M_\beta,
                 \mathcal H\otimes q_!p^*\mathcal K).
\end{align*}
Here $\mathsf P_{\alpha,\beta}$ is interpreted by the same
compact-support formula. Since $q^*\mathcal H\otimes p^*\mathcal K$
is bounded and mixed, Lemma~\ref{lem:convergent} shows that the
displayed complexes belong to $\Cq$.
Finally the projection formula for $x_L$, followed by pushforward
to $\pt$, gives
\[
 \RGC(\M_\alpha,x_{L!}\Q\otimes \mathcal K)
 \simeq x_L^*\mathcal K.\qedhere
\]
\end{proof}

\section{Construction of the character functor}
\label{sec:geometry}\label{sec:flags}

We construct the character functor $\Y_\ii$ by pairing mixed complexes on
representation stacks with Lusztig sheaves. The resulting complexes,
with prescribed shifts and Tate twists, form its graded components.
Lusztig's restriction formula gives each component of
$\Y_\ii(\mathcal F\star\mathcal G)$ a filtration whose layers are
tensor products of components of $\Y_\ii(\mathcal F)$ and
$\Y_\ii(\mathcal G)$, with explicit shifts and Tate twists. Taking Frobenius traces
then yields a character homomorphism to a quantum polynomial algebra.

Fix a word $\ii=(i_1,\ldots,i_m)$ in $Q_0$, that is,
a finite sequence of vertices of $Q$.
For $\mathbf a=(a_1,\ldots,a_m)\in\mathbb N^m$, set
\[
 A_k=a_k\eps_{i_k},\qquad
 \wt{\mathbf a}=\sum_k A_k.
\]
Repetitions in $\ii$ are allowed, and zero entries of $\mathbf a$ are
retained.

\subsection{Flag varieties and Lusztig sheaves}

Let $\mathfrak F_{\ii,\mathbf a}$ be the stack of representations equipped
with a descending flag of subrepresentations
\begin{equation}\label{geo:descending-flag}
 L=L^0\supseteq L^1\supseteq\cdots\supseteq L^m=0,
 \qquad \dim(L^{k-1}/L^k)=A_k.
\end{equation}
For a family over a scheme $S$, write $\mathcal V_i$ for the vector
bundle on $S$ at vertex $i$. The flag consists of subbundles
$\mathcal V_i^k\subseteq\mathcal V_i$ such that the arrow maps
$\mathcal V_h:\mathcal V_{s(h)}\to\mathcal V_{t(h)}$ satisfy
\[
 \mathcal V_h(\mathcal V_{s(h)}^k)\subseteq\mathcal V_{t(h)}^k
 \qquad(h\in Q_1,\ 0\leq k\leq m).
\]
The quotient $\mathcal V_i^{k-1}/\mathcal V_i^k$ is locally free of rank
$a_k$ if $i=i_k$, and is zero otherwise.  Since $Q$ has no loops, the $k$th
quotient is a family of copies of $S_{i_k}$.  The forgetful morphism
\[
 \pi_{\mathbf a}:\mathfrak F_{\ii,\mathbf a}
       \longrightarrow\M_{\wt{\mathbf a}}
\]
is representable and proper. Indeed, put $\alpha=\wt{\mathbf a}$
and use the quotient presentation $\M_\alpha=[E_\alpha/G_\alpha]$.
After base change along $E_\alpha\to\M_\alpha$, the morphism
$\pi_{\mathbf a}$ becomes the map $(x,L^\bullet)\mapsto x$ forgetting
an $x$-stable flag. This is Lusztig's proper flag projection
\cite[Lemma~1.6(d)]{Lusztig1991}.

Fix a representation $L=(L_i,L_h)$ of dimension $\alpha=\wt{\mathbf a}$,
with vector spaces $L_i$ at the vertices, and put
\[
 d_i^k=\sum_{\substack{l>k\\i_l=i}}a_l,
 \qquad i\in Q_0,\quad 0\leq k\leq m.
\]
Let $\operatorname{Fl}(d_i^\bullet;L_i)$ be the partial flag variety
parametrizing chains
\[
 L_i=L_i^0\supseteq L_i^1\supseteq\cdots\supseteq L_i^m=0,
 \qquad \dim L_i^k=d_i^k.
\]
Equal consecutive subspaces are allowed. For the representation point
$x_L:\pt\to\M_\alpha$, we have
\[
\begin{aligned}
 \Filt_{\ii,\mathbf a}(L)
 &=\pt\times_{x_L,\M_\alpha,\pi_{\mathbf a}}\mathfrak F_{\ii,\mathbf a}\\
 &=\left\{(L_i^\bullet)_i\in
       \prod_{i\in Q_0}\operatorname{Fl}(d_i^\bullet;L_i)
 \ \middle|\ 
 \begin{array}{l}
 L_h(L_{s(h)}^k)\subseteq L_{t(h)}^k\\
 \text{for all }h\in Q_1,\ 0\leq k\leq m
 \end{array}\right\}.
\end{aligned}
\]
The closed subscheme structure is cut out by the vanishing of the
composites
\[
 L_{s(h)}^k\hookrightarrow L_{s(h)}
 \xrightarrow{\,L_h\,}L_{t(h)}
 \twoheadrightarrow L_{t(h)}/L_{t(h)}^k
\]
for all $h\in Q_1$ and $0\leq k\leq m$.

Define
\begin{equation}\label{geo:test-complex}
 \mathcal K_{\mathbf a}=\pi_{\mathbf a!}\Q.
\end{equation}
We use the unshifted constant sheaf in this definition.  In particular,
$\mathcal K_{\mathbf a}$ is a bounded mixed complex. No assertion that it is
perverse is needed.

The Lusztig sheaves can also be obtained by iterated induction.
For $i\in Q_0$ and $a\in\mathbb N$, put
$\mathbf 1_{i,a}=\Q_{\M_{a\eps_i}}$.
When $a=0$, this is the induction unit $\mathbf1=\Q_{\M_0}$.
Using induction $\diamond$, we have
\begin{equation}\label{geo:iterated-tests}
 \mathcal K_{\mathbf a}\simeq
 \mathbf 1_{i_1,a_1}\diamond\cdots\diamond
 \mathbf 1_{i_m,a_m}.
\end{equation}
These are Lusztig sheaves in quotient-stack notation, with the unshifted
convention of \cite[\S2.2]{Lusztig1991}.

\subsection{Restriction of Lusztig sheaves}

For $\mathbf a,\mathbf b\in\mathbb N^m$, write
$A_k=a_k\eps_{i_k}$, $B_l=b_l\eps_{i_l}$, and set
\begin{equation}\label{geo:filtered-euler}
 \ef{\mathbf a}{\mathbf b}=\sum_{k\leq l}\eu{A_k}{B_l}.
\end{equation}

\begin{proposition}[Restriction of Lusztig sheaves]\label{geo:restriction-filtration}
Let $\wt{\mathbf c}=\alpha+\beta$. Lusztig's restriction formula gives
an isomorphism in $\Db_c((\M_\alpha\times\M_\beta)_{\bar k})$:
\begin{equation}\label{geo:restriction-decomposition}
 \bigl(\mathcal R_{\alpha,\beta}(\mathcal K_{\mathbf c})\bigr)_{\bar k}
 \simeq
 \bigoplus_{\substack{\mathbf a+\mathbf b=\mathbf c\\
                 \wt{\mathbf a}=\alpha,\ \wt{\mathbf b}=\beta}}
 (\mathcal K_{\mathbf a}\boxt\mathcal K_{\mathbf b})_{\bar k}
 [2\ef{\mathbf a}{\mathbf b}]\bigl(\ef{\mathbf a}{\mathbf b}\bigr).
\end{equation}
Over $k$, the same stratification gives a finite filtration of
$\mathcal R_{\alpha,\beta}(\mathcal K_{\mathbf c})$ in
$\Db_m(\M_\alpha\times\M_\beta)$, with layers
\begin{equation}\label{geo:restriction-layers}
 (\mathcal K_{\mathbf a}\boxt\mathcal K_{\mathbf b})
 [2\ef{\mathbf a}{\mathbf b}]\bigl(\ef{\mathbf a}{\mathbf b}\bigr),
\end{equation}
indexed by the same pairs $(\mathbf a,\mathbf b)$.
\end{proposition}

\begin{proof}
We compute the shifts and Tate twists in
\eqref{geo:restriction-layers} under our quotient-stack conventions.
The remaining steps follow Lusztig's restriction argument
\cite[\S\S4.3--4.7]{Lusztig1991}, expressed in quotient-stack notation.

In the Hall correspondence, an extension
$0\to N\to L\to M\to0$, with $\dim M=\alpha$ and $\dim N=\beta$,
together with a flag $L^\bullet$ of type $\mathbf c$ induces the flags
\[
 M^r=(L^r+N)/N,\qquad N^r=N\cap L^r.
\]
The types $\mathbf a$ of $M^\bullet$ and $\mathbf b$ of $N^\bullet$
index Lusztig's locally closed strata \cite[\S4.3]{Lusztig1991}.
Fix the underlying graded spaces $N\subset L$ and $M=L/N$.
Write $\widetilde E_{\mathbf a}$, $\widetilde E_{\mathbf b}$,
and $\widetilde E_{\mathbf c}$ for the varieties of representations
with stable flags of the indicated types on $M$, $N$, and $L$,
respectively. Let $Z^{\mathrm{var}}_{\mathbf a,\mathbf b}
\subset\widetilde E_{\mathbf c}$ be the locus where $N$ is stable
and the induced flags on $M,N$ have types $\mathbf a,\mathbf b$.
The morphism
\[
 \begin{aligned}
 f_{\mathbf a,\mathbf b}:Z^{\mathrm{var}}_{\mathbf a,\mathbf b}
 &\longrightarrow\widetilde E_{\mathbf a}\times\widetilde E_{\mathbf b},\\
 (x,L^\bullet)&\longmapsto
 \bigl((x_{L/N},(L^\bullet+N)/N),
       (x|_N,N\cap L^\bullet)\bigr)
 \end{aligned}
\]
is a vector bundle of rank
\[
 d(\mathbf a,\mathbf b)=
 \sum_{k\leq l}\sum_{h\in Q_1}(A_k)_{s(h)}(B_l)_{t(h)}
 +\sum_{k>l}\sum_{i\in Q_0}(A_k)_i(B_l)_i
\]
by \cite[Lemma~4.4]{Lusztig1991}.
Here $(A_k)_i=a_k\delta_{i,i_k}$ and
$(B_l)_i=b_l\delta_{i,i_l}$ are components of dimension vectors.
Passing to the extension stack involves quotienting by the
unipotent group
\[
 U=\prod_{i\in Q_0}\operatorname{Hom}(M_i,N_i),
 \qquad u=\dim U=\sum_i\alpha_i\beta_i,
\]
where $M_i,N_i$ are the vertex spaces of $M,N$.
Let $\mathfrak Z_{\mathbf a,\mathbf b}$
denote the stack of these flagged extensions whose induced flags
have types $\mathbf a,\mathbf b$. The induced-flag morphism
\[
 h_{\mathbf a,\mathbf b}:\mathfrak Z_{\mathbf a,\mathbf b}
 \longrightarrow
 \mathfrak F_{\ii,\mathbf a}\times\mathfrak F_{\ii,\mathbf b}
\]
has compact-support pushforward
\[
 h_{\mathbf a,\mathbf b!}\Q
 \simeq\Q[-2(d(\mathbf a,\mathbf b)-u)]
              (-(d(\mathbf a,\mathbf b)-u)).
\]
Compactly supported pushforward along the rank-$d$ vector bundle
$f_{\mathbf a,\mathbf b}$ contributes the shift and twist $[-2d](-d)$.
The quotient by $U$ contributes $[2u](u)$.
Since $\alpha=\sum_k A_k$ and $\beta=\sum_l B_l$, we have
\[
 u=\sum_i\alpha_i\beta_i
   =\sum_{k,l}\sum_{i\in Q_0}(A_k)_i(B_l)_i.
\]
Subtracting $u$ cancels the terms with $k>l$ in
$d(\mathbf a,\mathbf b)$. By the definition of the Euler form,
\[
 \begin{aligned}
 d(\mathbf a,\mathbf b)-u
 &=\sum_{k\leq l}\left(
       \sum_{h\in Q_1}(A_k)_{s(h)}(B_l)_{t(h)}
       -\sum_{i\in Q_0}(A_k)_i(B_l)_i\right)\\
 &=-\sum_{k\leq l}\eu{A_k}{B_l}
  =-\ef{\mathbf a}{\mathbf b}.
 \end{aligned}
\]
This calculation gives the shifts and Tate twists in
\eqref{geo:restriction-layers}. With these conventions, the mixed
filtration is the one constructed in \cite[\S4.6]{Lusztig1991}.
The flag stacks are smooth and their projections are proper, so
the layers in \eqref{geo:restriction-layers}, and hence the filtration
stages, are pure of weight zero.
Since $\M_\alpha\times\M_\beta$ has affine stabilizers,
\cite[Corollary~3.10(iii)]{SunDecomposition2012} implies that
these triangles split over $\bar k$. This gives
\eqref{geo:restriction-decomposition}, the quotient-stack form of
\cite[Proposition~4.2(b) and Lemma~4.7]{Lusztig1991}.
\end{proof}

\begin{remark}
More generally, Lusztig's restriction functor preserves geometric
semisimplicity for equivariant mixed complexes
\cite[Proposition~2.10]{FangLanXiao2023}.
\end{remark}

\subsection{Normalization associated with flag types}

For a flag type $(\ii,\mathbf a)$ of weight $\alpha$, put
\begin{equation}\label{char:normalizing-exponents}
 R(\mathbf a)=\sum_{k<l}\eu{A_l}{A_k}.
\end{equation}
The following dimension calculation explains this normalization.

\begin{lemma}\label{geo:n-geometric}
The stacks $\mathfrak F_{\ii,\mathbf a}$ and $\M_\alpha$ are smooth, and
\[
 R(\mathbf a)=\dim\mathfrak F_{\ii,\mathbf a}-\dim\M_\alpha.
\]
\end{lemma}

\begin{proof}
Fix a graded vector space $W$ of dimension $\alpha$ and a flag of
type $(\ii,\mathbf a)$. Let $P\subset G_\alpha$ be its stabilizer
and $A$ the space of arrow maps preserving the flag. Then
$\mathfrak F_{\ii,\mathbf a}\simeq[A/P]$, where $A$ is an affine
space and $P$ is smooth. Moreover,
\[
 \begin{aligned}
 \dim A&=\sum_{k<l}\sum_{h\in Q_1}
                (A_k)_{s(h)}(A_l)_{t(h)},\\
 \dim P&=\sum_k a_k^2+\sum_{k<l}\sum_{i\in Q_0}(A_k)_i(A_l)_i.
 \end{aligned}
\]
Thus
\[
 \dim\mathfrak F_{\ii,\mathbf a}
 =-\sum_k a_k^2-\sum_{k<l}\eu{A_k}{A_l}.
\]
Since $\M_\alpha=[E_\alpha/G_\alpha]$ is smooth of dimension
$-\eu{\alpha}{\alpha}$, the assertion follows.
\end{proof}

For $\mathbf c=\mathbf a+\mathbf b$, combine the induction shift
in \eqref{eq:star} with the normalization $\bal{R(\mathbf c)}$.
Relative to the separate normalizations $\bal{R(\mathbf a)}$ and
$\bal{R(\mathbf b)}$, the remaining shift in
\eqref{geo:restriction-layers} is described by the skew-symmetric form
\begin{equation}\label{geo:ordered-form}
 L_\ii(\mathbf a,\mathbf b)
 =\sum_{k<l}\bigl(\eu{A_l}{B_k}-\eu{B_l}{A_k}\bigr).
\end{equation}

\begin{lemma}\label{geo:normalization-identity}
For $\wt{\mathbf a}=\alpha$ and $\wt{\mathbf b}=\beta$,
\begin{equation}\label{geo:normalization-equation}
 2\ef{\mathbf a}{\mathbf b}-2\eu{\alpha}{\beta}
 +R(\mathbf a+\mathbf b)-R(\mathbf a)-R(\mathbf b)
 =-L_\ii(\mathbf a,\mathbf b).
\end{equation}
\end{lemma}

\begin{proof}
By bilinearity,
\[
 R(\mathbf a+\mathbf b)-R(\mathbf a)-R(\mathbf b)
 =\sum_{k<l}\bigl(\eu{A_l}{B_k}+\eu{B_l}{A_k}\bigr).
\]
On the other hand,
\[
 \ef{\mathbf a}{\mathbf b}-\eu{\alpha}{\beta}
 =-\sum_{k<l}\eu{A_l}{B_k}.
\]
Substitution gives the assertion.
\end{proof}

\subsection{Character functor and multiplicativity}

\begin{definition}[Character functor]\label{geo:character-definition}

Let $\mathcal T_{\ii}$ be the category of finitely supported
$\mathbb N^m$-graded objects in $\Cq$, with degree-preserving morphisms.
Shifts and distinguished triangles are defined componentwise.
The character functor is
\[
 \begin{aligned}
 \Y_{\ii}:\mathscr D_Q&\longrightarrow\mathcal T_{\ii},\\
 \mathcal F&\longmapsto
       \bigl(\Y_{\ii}(\mathcal F)_{\mathbf a}\bigr)_{\mathbf a\in\mathbb N^m}.
 \end{aligned}
\]
For $\mathcal F\in\Db_m(\M_\alpha)$ its components are defined by

\begin{equation}\label{geo:character}
 \Y_{\ii}({\mathcal F})_{\mathbf a}
 =\begin{cases}
   \RGC(\M_\alpha,{\mathcal F}\otimes \mathcal K_{\mathbf a})
                           \bal{R(\mathbf a)},
       &\wt{\mathbf a}=\alpha,\\
   0,&\wt{\mathbf a}\ne\alpha.
  \end{cases}
\end{equation}
The definition extends to all of $\mathscr D_Q$ by finite direct sums.
\end{definition}

There are only finitely many $\mathbf a$ of a fixed weight $\alpha$.
Moreover, ${\mathcal F}\otimes \mathcal K_{\mathbf a}$ is bounded and mixed on the
finite-type quotient stack $\M_\alpha$.
Lemma~\ref{lem:convergent} therefore applies to each component of
\eqref{geo:character}.

Tensoring with a fixed complex and compact-support pushforward are
exact, so $\Y_\ii$ is an exact functor.
We equip $\mathcal T_\ii$ with the product

\begin{equation}\label{geo:target-product}
 (U\mathbin{\odot}W)_{\mathbf c}
 =\bigoplus_{\mathbf a+\mathbf b=\mathbf c}
    (U_{\mathbf a}\otimes W_{\mathbf b})
       \bal{-L_\ii(\mathbf a,\mathbf b)}.
\end{equation}

This product is biexact and preserves finite support, with unit $\Q$
in grading degree zero. With the associativity constraints induced
by the normalization in Lemma~\ref{geo:normalization-identity},
$(\mathcal T_\ii,\odot)$ is a monoidal category. The induced associative unital multiplication on
$K_0(\mathcal T_\ii)$ is $[U][W]=[U\odot W]$.
The Frobenius-trace realization of this multiplication is described in
Section~\ref{geo:realizations}.

\begin{theorem}[Multiplicative filtration]\label{geo:main-theorem}
For ${\mathcal F}\in\Db_m(\M_\alpha)$, ${\mathcal G}\in\Db_m(\M_\beta)$, and
$\wt{\mathbf c}=\alpha+\beta$, the complex
$\Y_{\ii}({\mathcal F}\star {\mathcal G})_{\mathbf c}$ admits a finite filtration,
functorial in ${\mathcal F}$ and ${\mathcal G}$, with layers
\begin{equation}\label{geo:main-layers}
 \bigl(\Y_{\ii}({\mathcal F})_{\mathbf a}\otimes
              \Y_{\ii}({\mathcal G})_{\mathbf b}\bigr)
      \bal{-L_\ii(\mathbf a,\mathbf b)},
\end{equation}
indexed by $\mathbf a+\mathbf b=\mathbf c$ with
$\wt{\mathbf a}=\alpha$ and $\wt{\mathbf b}=\beta$.
Thus
\[
 \operatorname{gr}\!\left(
   \Y_{\ii}(\mathcal F\star\mathcal G)_{\mathbf c}\right)
 \simeq
 \bigl(\Y_{\ii}(\mathcal F)\odot
       \Y_{\ii}(\mathcal G)\bigr)_{\mathbf c},
\]
where $\operatorname{gr}$ denotes the direct sum of the layers of this filtration.
\end{theorem}

\begin{proof}
Set $\mathcal H=\mathcal F\boxt\mathcal G$.
By Lemma~\ref{geo:pairing-transposition},
\[
 \begin{aligned}
 \RGC\bigl(\M_{\alpha+\beta},
   (\mathcal F\diamond\mathcal G)\otimes\mathcal K_{\mathbf c}\bigr)
 &\simeq \mathsf P_{\alpha,\beta}
   (\mathcal H,\mathcal R_{\alpha,\beta}\mathcal K_{\mathbf c})\\
 &=\RGC\bigl(\M_\alpha\times\M_\beta,
   \mathcal H\otimes\mathcal R_{\alpha,\beta}\mathcal K_{\mathbf c}\bigr).
 \end{aligned}
\]
Apply $\RGC(\M_\alpha\times\M_\beta,\mathcal H\otimes-)$ to
the filtration in Proposition~\ref{geo:restriction-filtration} and include
the shifts defining $\star$ and $\Y_{\ii}(-)_{\mathbf c}$.
Writing
\[
 d_{\mathbf a,\mathbf b}
 =2\ef{\mathbf a}{\mathbf b}-2\eu{\alpha}{\beta}+R(\mathbf c),
\]
the resulting $(\mathbf a,\mathbf b)$ layer of
$\Y_{\ii}(\mathcal F\star\mathcal G)_{\mathbf c}$ is
\[
 \begin{aligned}
 &\RGC\bigl(\M_\alpha\times\M_\beta,
   \mathcal H\otimes(\mathcal K_{\mathbf a}\boxt\mathcal K_{\mathbf b})\bigr)
   \bal{d_{\mathbf a,\mathbf b}}\\
 &\qquad\simeq
 \bigl(\RGC(\M_\alpha,{\mathcal F}\otimes \mathcal K_{\mathbf a})
       \otimes\RGC(\M_\beta,{\mathcal G}\otimes \mathcal K_{\mathbf b})\bigr)
       \bal{d_{\mathbf a,\mathbf b}}\\
 &\qquad\simeq
 \bigl(\Y_{\ii}(\mathcal F)_{\mathbf a}\otimes
       \Y_{\ii}(\mathcal G)_{\mathbf b}\bigr)
       \bal{d_{\mathbf a,\mathbf b}-R(\mathbf a)-R(\mathbf b)}\\
 &\qquad\simeq
 \bigl(\Y_{\ii}(\mathcal F)_{\mathbf a}\otimes
       \Y_{\ii}(\mathcal G)_{\mathbf b}\bigr)
       \bal{-L_\ii(\mathbf a,\mathbf b)}.
 \end{aligned}
\]
The last isomorphism uses $\mathbf c=\mathbf a+\mathbf b$ and
Lemma~\ref{geo:normalization-identity}.
The functors used above are exact and natural in $\mathcal F,\mathcal G$,
so the resulting filtration is functorial in $\mathcal F$ and $\mathcal G$.
\end{proof}

\subsection{Quantum tori and normalized monomials}
\label{sec:quantum-tori}

Fix a word $\ii=(i_1,\ldots,i_m)$ in $Q_0$. Use $\mathbf e_k$ for the
standard basis of $\mathbb Z^m$. For an exponent vector
$\mathbf a=(a_1,\ldots,a_m)\in\mathbb Z^m$, put
\[
 \wt{\mathbf a}=\sum_{k=1}^m a_k\eps_{i_k}\in\mathbb Z^{Q_0}.
\]
Thus $\mathbf e_k$ and
$\eps_{i_k}$ belong to the word and vertex lattices, respectively.
Let $\mathcal P_{\ii}$ be the quantum polynomial algebra over $\mathbb C$ generated by
$t_1,\ldots,t_m$, subject to

\begin{equation}\label{eq:ordered-relations}
 t_l t_k=w^{2\eu{\eps_{i_l}}{\eps_{i_k}}}t_k t_l
 \qquad(k<l).
\end{equation}

For $\mathbf a\in\mathbb N^m$, define the ordered monomial
$t_{\mathrm{ord}}^{\mathbf a}=t_1^{a_1}\cdots t_m^{a_m}$, with the
factors written in the order $t_1,\ldots,t_m$.
For $\mathbf a,\mathbf b\in\mathbb N^m$, we have

\begin{equation}\label{eq:ordered-product}
 t_{\mathrm{ord}}^{\mathbf a}t_{\mathrm{ord}}^{\mathbf b}
 =w^{2\sum_{k<l}\eu{A_l}{B_k}}
    t_{\mathrm{ord}}^{\mathbf a+\mathbf b}.
\end{equation}

Inverting the generators gives the quantum torus $\TT_{\ii}$, in which
the ordered monomials and their multiplication formula extend to
$\mathbf a,\mathbf b\in\mathbb Z^m$.

Following the based-monomial convention of \cite[\S4.1]{BZ2005}, set
\begin{equation}\label{eq:based-monomial}
 t^{\mathbf a}=w^{R(\mathbf a)}t_{\mathrm{ord}}^{\mathbf a}.
\end{equation}
For the integral skew-symmetric form $L_\ii$ in
\eqref{geo:ordered-form}, one has
$L_\ii(\mathbf e_k,\mathbf e_l)=-\eu{\eps_{i_l}}{\eps_{i_k}}$ for $k<l$.
Then $t^{\mathbf a}t^{\mathbf b}
=w^{L_\ii(\mathbf a,\mathbf b)}t^{\mathbf a+\mathbf b}$.

\subsection{Character homomorphisms and Frobenius traces}
\label{geo:realizations}

The triangle relations turn the filtration of
Theorem~\ref{geo:main-theorem} into a multiplicative identity in
the triangulated Grothendieck groups.

\begin{corollary}\label{geo:K0-homomorphism}
The exact functor $\Y_{\ii}$ induces a unital ring homomorphism
\[
 K_0\!\left(\bigoplus_\alpha\Db_m(\M_\alpha)\right)
       \longrightarrow K_0(\mathcal T_{\ii}).
\]
\end{corollary}

\begin{proof}
Exactness gives a homomorphism of abelian groups.  Applying the
triangle relations to the finite filtration in
Theorem~\ref{geo:main-theorem} proves compatibility with products.
The zero representation has a unique flag of type zero, so the Hall
unit maps to $\Q$ in degree zero.
\end{proof}

For $U\in\Cq$, put
\[
 \Tr(U)=\sum_{j\in\mathbb Z}(-1)^j
   \iota\!\left(\Tr(\Fr;H^j(U_{\bar k}))\right)\in\mathbb C.
\]
The sum converges absolutely by the definition of $\Cq$.
This trace is additive on distinguished triangles and
multiplicative under tensor product.  It satisfies
\[
 \Tr(U\bal{d})=(-1)^dq^{-d/2}\Tr(U)
                    =w^{-d}\Tr(U),\qquad d\in\mathbb Z,
\]
with the square-root convention of \eqref{eq:paired-shift}.
In $\mathcal P_{\ii}$, the monomials $t^{\mathbf a}$ multiply with the
factor $w^{L_\ii(\mathbf a,\mathbf b)}$. This is the trace factor of the
paired shift in \eqref{geo:target-product}. Define
\begin{equation}\label{geo:trace-character}
 \Psi_{\ii}([{\mathcal F}])
    =\sum_{\mathbf a}\Tr\bigl(\Y_{\ii}({\mathcal F})_{\mathbf a}\bigr)
                   t^{\mathbf a}.
\end{equation}

\begin{corollary}\label{geo:trace-homomorphism}
The map $\Psi_{\ii}:K_0(\mathscr D_Q)\to\mathcal P_{\ii}$ is a unital
ring homomorphism. It extends to a $\mathbb C$-algebra homomorphism
\[
 \begin{aligned}
 \Psi_{\ii,\mathbb C}:K_0(\mathscr D_Q)\otimes_{\mathbb Z}\mathbb C
   &\longrightarrow\mathcal P_{\ii},\\
 x\otimes z&\longmapsto z\Psi_{\ii}(x).
 \end{aligned}
\]
In particular,
\begin{equation}\label{geo:trace-multiplication}
 \Psi_{\ii}([{\mathcal F}\star {\mathcal G}])
        =\Psi_{\ii}([{\mathcal F}])\Psi_{\ii}([{\mathcal G}]).
\end{equation}
\end{corollary}

\begin{proof}
The graded trace is additive on triangles, and the calculation above
shows that it respects multiplication.
It therefore defines a ring homomorphism
$K_0(\mathcal T_{\ii})\to\mathcal P_{\ii}$.  Composing with
Corollary~\ref{geo:K0-homomorphism} proves the assertion.
\end{proof}

Thus taking Frobenius traces of the filtration in
Theorem~\ref{geo:main-theorem} gives the multiplicative identity
\eqref{geo:trace-multiplication}, yielding the character homomorphism
$\Psi_{\ii}$ with values in the quantum polynomial algebra
$\mathcal P_{\ii}\subset\TT_{\ii}$.

\section{Comparison with character maps}
\label{sec:characters}

We compare the homomorphism on the dual Hall algebra induced by the character functor of
Section~\ref{sec:flags} with the character homomorphism of
Berenstein and Rupel. We then recover Reineke's integration map, the Feigin
homomorphism, and quantum cluster characters from this construction.

\subsection{Trace functions and Hall algebras}

All numerical algebras below are over $\mathbb C$.
For representations $M,N,L$, set $a_L=|\Aut(L)(k)|$ and
\[
 F^L_{MN}=\#\{U\leq L\mid U\simeq N,\ L/U\simeq M\},
\]
where $U\leq L$ denotes a subrepresentation. The untwisted Hall algebra has
basis $u_L$ and product
\begin{equation}\label{eq:raw-hall}
 u_Mu_N=\sum_{[L]}F^L_{MN}u_L.
\end{equation}
The Hall algebra $\HH_w$ associated with $\star$ has basis $[L]$ and product
\[
 [M][N]=w^{2\eu{\alpha}{\beta}}\sum_{[L]}F^L_{MN}[L],
 \qquad \alpha=\dim M,\quad\beta=\dim N.
\]
These are the respective structure constants of the trace functions of
$\diamond$ and $\star$. Write
$\HH_w^*=\bigoplus_\alpha(\HH_w)_\alpha^*$ for the graded
dual, with basis $\delta_L([W])=\delta_{[L],[W]}$.

Let $\mathrm{CF}_\alpha$ be the space of $\mathbb C$-valued functions
on the $k$-isomorphism classes of representations of dimension $\alpha$,
and set
\[
 \mathrm{CF}(Q)=\bigoplus_{\alpha\in\mathbb N^{Q_0}}\mathrm{CF}_\alpha.
\]
Thus a function in $\mathrm{CF}(Q)$ is supported in finitely many
dimension vectors. Each dimension vector has only finitely many
$k$-isomorphism classes, so all sums below are finite.
The stack measure in \eqref{eq:stack-trace} defines a linear isomorphism
\[
 \begin{aligned}
 J:\mathrm{CF}(Q)&\xrightarrow{\sim}\HH_w^*,\\
 J(f)\!\left(\sum_{[L]}h_L[L]\right)
   &=\sum_{[L]}\frac{f(L)h_L}{a_L}.
 \end{aligned}
\]
Equivalently,
\begin{equation}\label{eq:dual-pairing}
 J(f)=\sum_{[L]}\frac{f(L)}{a_L}\delta_L,
 \qquad J^{-1}(\delta_L)=a_L\mathbf1_L,
\end{equation}
where $\mathbf1_L$ is the characteristic function of the isomorphism
class of $L$.
Define the product on $\HH_w^*$ by $J(f)J(g)=J(f\star g)$. Then
\begin{equation}\label{eq:dual-product}
 \delta_M\delta_N
 =w^{2\eu{\alpha}{\beta}}
   \sum_{[L]}F^L_{MN}\frac{a_Ma_N}{a_L}\delta_L.
\end{equation}
The factor $w^{2\eu{\alpha}{\beta}}=q^{\eu{\alpha}{\beta}}$
comes from the shift in $q^!$.

\subsection{The induced Hall algebra homomorphism}

Let $L$ be a representation over $k$ of dimension $\alpha$, and write
\[
 x_L\colon\pt\longrightarrow\M_\alpha,
 \qquad A_L=(x_L)_!\Q.
\]
\begin{proposition}\label{char:point-fibre}
For every flag type $\mathbf a$ of weight $\alpha$, there is a natural
isomorphism of Weil complexes
\begin{equation}\label{char:point-complex}
 \RGC(\M_\alpha,A_L\otimes \mathcal K_{\mathbf a})
 \simeq
 \RGC(\Filt_{\ii,\mathbf a}(L),\Q).
\end{equation}
Moreover, the trace function of $A_L$ is $a_L\mathbf1_L$.
\end{proposition}

\begin{proof}
The square
\[
 \begin{tikzcd}
 \Filt_{\ii,\mathbf a}(L) \arrow[r] \arrow[d]
   & \mathfrak F_{\ii,\mathbf a} \arrow[d,"\pi_{\mathbf a}"] \\
 \pt \arrow[r,"x_L"'] & \M_\alpha
 \end{tikzcd}
\]
is cartesian: fixing an identification of the underlying representation
with $L$ identifies the fibre with flags of subrepresentations in $L$.
The projection formula for $x_L$ and proper base change give
\[
 \RGC(\M_\alpha,A_L\otimes\mathcal K_{\mathbf a})
 \simeq x_L^*\mathcal K_{\mathbf a}
 \simeq\RGC(\Filt_{\ii,\mathbf a}(L),\Q),
\]
proving \eqref{char:point-complex}.
For a representation $W$, base change also gives
\[
 x_W^*A_L\simeq\RGC(\operatorname{Isom}_Q(L,W),\Q).
\]
The trace formula therefore yields
\[
 f_{A_L}(W)=\Tr(\Fr,x_W^*A_L)
 =\#\operatorname{Isom}_Q(L,W)(k)
 =\begin{cases}
 a_L,&W\simeq L,\\
 0,&W\not\simeq L,
 \end{cases}
\]
so $f_{A_L}=a_L\mathbf1_L$.
\end{proof}

Thus the coefficients of $\Y_\ii(A_L)$ are the actual complexes
\begin{equation}\label{char:Y-at-V}
 \Y_\ii(A_L)_{\mathbf a}
 \simeq
 \RGC(\Filt_{\ii,\mathbf a}(L),\Q)
 \bal{R(\mathbf a)}.
\end{equation}
In particular, no division by $a_L$ occurs in their traces.

By \eqref{eq:dual-pairing}, the trace function
$a_L\mathbf1_L$ of $A_L$ maps to $\delta_L$ under $J$.
The point object therefore differs from the constant complex on
its orbit, whose trace is $\mathbf1_L$.

We now identify the Hall algebra homomorphism obtained by taking
Frobenius traces of the multiplicative filtration of $\Y_\ii$.
For $\mathcal F\in\Db_m(\M_\alpha)$, put
$f_{\mathcal F}(L)=\Tr(\Fr,x_L^*\mathcal F)$.

\begin{proposition}\label{char:Hall-realization}
The homomorphism in \eqref{geo:trace-character} factors through
$[\mathcal F]\mapsto J(f_{\mathcal F})$, inducing a unital algebra
homomorphism, also denoted by
\[
 \Psi_\ii:\HH_w^*\longrightarrow\mathcal P_\ii\subseteq\TT_\ii.
\]
It satisfies
\[
 \Psi_\ii\bigl(J(f_{\mathcal F})\bigr)
 =\sum_{\mathbf a}\Tr\bigl(\Y_\ii(\mathcal F)_{\mathbf a}\bigr)
 t^{\mathbf a}.
\]
\end{proposition}

\begin{proof}
For $\wt{\mathbf a}=\alpha$, the definition of $\Y_\ii$,
the stack trace formula \eqref{eq:stack-trace}, and proper base change give
\begin{equation}\label{char:trace-factorization}
 \begin{aligned}
 \Tr\bigl(\Y_\ii(\mathcal F)_{\mathbf a}\bigr)
 &=w^{-R(\mathbf a)}\Tr\bigl(\Fr,
     \RGC(\M_\alpha,\mathcal F\otimes\mathcal K_{\mathbf a})\bigr)\\
 &=w^{-R(\mathbf a)}\sum_{[L]}\frac{f_{\mathcal F}(L)}{a_L}
     \Tr\bigl(\Fr,x_L^*\mathcal K_{\mathbf a}\bigr)\\
 &=w^{-R(\mathbf a)}\sum_{[L]}\frac{f_{\mathcal F}(L)}{a_L}
     |\Filt_{\ii,\mathbf a}(L)(k)|.
 \end{aligned}
\end{equation}
Since $J(f_{\mathcal F})=\sum_{[L]}f_{\mathcal F}(L)\delta_L/a_L$,
the result depends only on $J(f_{\mathcal F})$ and defines the asserted
linear map. Since $J$ respects convolution and
$J(f_{A_L})=\delta_L$ for every $L$, multiplicativity and the unit
follow from Corollary~\ref{geo:trace-homomorphism}.
\end{proof}

\subsection{Character and integration maps}

\leavevmode In the basis $t^{\mathbf a}$, the coefficients of $\Psi_\ii(\delta_L)$
agree with those of the quantum cluster $\ii$-character of Berenstein and Rupel
\cite[\S2]{BR2015} when their square root of $q$ is taken to be $w=-\sqrt q$.

\begin{corollary}[Quantum cluster $\ii$-characters]\label{char:BR}
For a representation $L$ of dimension $\alpha$,
\begin{align}
 X_{L,\ii}:=\Psi_\ii(\delta_L)
 &=\sum_{\wt{\mathbf a}=\alpha}
 \Tr\bigl(\Y_\ii(A_L)_{\mathbf a}\bigr)t^{\mathbf a}
 \label{char:Z-trace}\\
 &=\sum_{\wt{\mathbf a}=\alpha}
 w^{-R(\mathbf a)}|\Filt_{\ii,\mathbf a}(L)(k)|\,t^{\mathbf a}
 \label{char:BR-based}\\
 &=\sum_{\wt{\mathbf a}=\alpha}
 |\Filt_{\ii,\mathbf a}(L)(k)|\,t_1^{a_1}\cdots t_m^{a_m}.
 \label{char:BR-ordered}
\end{align}
\end{corollary}

\begin{proof}
Take Frobenius traces in \eqref{char:Y-at-V} and use
\eqref{eq:based-monomial}.
\end{proof}

\begin{remark}\label{char:comparison}
Fix a square root $w^{1/2}\in\mathbb C$. On the graded vector spaces
$\HH_w^*$ and $\mathcal P_\ii$, define
\[
 \begin{aligned}
 \delta_M\cdot_{\rm BR}\delta_N
 &=w^{\frac12(\eu{\beta}{\alpha}-\eu{\alpha}{\beta})}
       (\delta_M\delta_N),\\
 t^{\mathbf a}\cdot_{\rm BR}t^{\mathbf b}
 &=w^{L_\ii(\mathbf a,\mathbf b)
       +\frac12(\eu{\beta}{\alpha}-\eu{\alpha}{\beta})}
       t^{\mathbf a+\mathbf b},
 \end{aligned}
\]
where $\alpha=\dim M=\wt{\mathbf a}$ and
$\beta=\dim N=\wt{\mathbf b}$. In the first formula,
$\delta_M\delta_N$ is the Hall product \eqref{eq:dual-product}.
The map $\Psi_\ii:(\HH_w^*,\cdot_{\rm BR})
\to(\mathcal P_\ii,\cdot_{\rm BR})$
remains an algebra homomorphism. It agrees with the character
homomorphism of \cite{BR2015} under the identification of normalized
monomials and
\[
 \delta_L\longmapsto[L]^*
 =w^{(\ellone{\alpha}-\eu{\alpha}{\alpha})/2}\delta_L^{\rm BR},
 \qquad\alpha=\dim L.
\]
Here $\delta_L^{\rm BR}$ is the basis element dual to $[L]$
in their dual Hall algebra.
\end{remark}

An ordering $(i_1,\ldots,i_n)$ of the vertices of $Q$ is
\emph{source-adapted} if every arrow $i_r\to i_s$ satisfies $r<s$.
Equivalently, after deleting $i_1,\ldots,i_{r-1}$, no arrow enters $i_r$
in the remaining quiver. Since $Q$ is acyclic, we may choose such an
ordering and relabel the vertices so that $\iiO=(1,\ldots,n)$.

\begin{corollary}[Integration for a source-adapted word]\label{char:source}
For every $\alpha$, the morphism
\[
 \pi_\alpha:\mathfrak F_{\iiO,\alpha}\longrightarrow\M_\alpha,
 \qquad (L,L^\bullet)\longmapsto L,
\]
is an isomorphism. Hence $\mathcal K_\alpha\simeq\Q$ and,
for $\mathcal F\in\Db_m(\M_\alpha)$,
\[
 \RGC(\M_\alpha,\mathcal F\otimes\mathcal K_\alpha)
 \simeq\RGC(\M_\alpha,\mathcal F).
\]
For a representation $L$ of dimension $\alpha$,
\begin{equation}\label{char:source-basis}
 \Psi_{\iiO}(\delta_L)=t^\alpha,
 \qquad X_{L,\iiO}=t^\alpha.
\end{equation}
\end{corollary}

\begin{proof}
A flag of type $(\iiO,\alpha)$ must satisfy
\[
 L^{>j}_r=\begin{cases}0&r\le j,\\L_r&r>j.\end{cases}
\]
Since every arrow $h$ satisfies $s(h)<t(h)$, each $L^{>j}$ is a
subrepresentation. Thus $L\mapsto(L,L^{>\bullet})$ is the inverse of
$\pi_\alpha$, and the same construction applies to families. Moreover,
\[
 R(\alpha)=\sum_{k<l}\alpha_k\alpha_l\eu{\eps_l}{\eps_k}=0.
\]
Now \eqref{char:source-basis} follows from
Corollary~\ref{char:BR}.
\end{proof}

\begin{corollary}[Reineke's integration map]\label{rei:trace}
Let $\mathbb A_Q$ be the $\mathbb C$-algebra defined by
\[
 \mathbb A_Q=\bigoplus_{\alpha\in\mathbb N^{Q_0}}\mathbb C z^\alpha,
 \qquad z^\alpha z^\beta=q^{-\eu{\alpha}{\beta}}z^{\alpha+\beta}.
\]
The linear map from the untwisted Hall algebra
\eqref{eq:raw-hall} to $\mathbb A_Q$ defined by
\begin{equation}\label{eq:reineke-trace}
 \int_Q u_L=\frac{z^{\dim L}}{a_L}
\end{equation}
is an algebra homomorphism \cite[Lemma~3.3]{Reineke2006}.
\end{corollary}

\begin{proof}
For this word, the multiplicative filtration of
Theorem~\ref{geo:main-theorem} has a single layer. Undoing the
convolution shift and taking
Frobenius traces gives
\[
 \sum_{[L]}\frac{F^L_{MN}}{a_L}
 =\frac{q^{-\eu{\alpha}{\beta}}}{a_Ma_N},
 \qquad \alpha=\dim M,\quad\beta=\dim N.
\]
This is the multiplicativity identity for \eqref{eq:reineke-trace}.
\end{proof}

The counting identity above is the sum over middle terms of the
Riedtmann--Peng formula \cite[\S II.4]{Ringel1996}:
\[
 \frac{a_Ma_N}{a_L}F^L_{MN}
 =\frac{|\Ext^1_Q(M,N)_L|}{|\Hom_Q(M,N)|},
\]
where $\Ext^1_Q(M,N)_L$ consists of extension classes with middle
term isomorphic to $L$. Summing over $[L]$ gives
\[
 \sum_{[L]}\frac{F^L_{MN}}{a_L}
 =\frac{|\Ext^1_Q(M,N)|}{a_Ma_N|\Hom_Q(M,N)|}
 =\frac{q^{-\eu{\alpha}{\beta}}}{a_Ma_N}.
\]

\begin{remark}[The Feigin homomorphism]\label{char:Feigin}
For any word $\ii$, a flag of a simple representation $S_j$ is a point
of type $\mathbf e_k$ with $i_k=j$, and its normalization is trivial.
Thus
\begin{equation}\label{char:Feigin-generators}
 \Psi_\ii(\delta_{S_j})=\sum_{k:i_k=j}t_k.
\end{equation}
With the source and target products of Remark~\ref{char:comparison},
the restriction to the subalgebra generated by these simple classes
is the homomorphism introduced by Feigin. For its Hall algebra interpretation, see
\cite{Reineke2001}.
\end{remark}

\subsection{Quantum cluster characters}
\label{sec:cluster}

For this ordering $\iiO$, put
$\ii=(\iiO,\iiO)$, following \cite[\S5]{BR2015}.

\begin{proposition}[Berenstein--Rupel {\cite[Theorem~5.2]{BR2015}}]
For every representation $L$ of dimension $\alpha$ and every
$0\le e\le\alpha$, the map
\begin{equation}\label{cc:flag-gr-isomorphism}
 \Filt_{\ii,(\alpha-e,e)}(L)\xrightarrow{\sim}\Gr_e(L),
 \qquad L^\bullet\longmapsto L^n
\end{equation}
is a bijection. Here $L^n$ is the $n$th term of the flag
$L=L^0\supseteq\cdots\supseteq L^{2n}=0$, and $\dim L^n=e$.
\end{proposition}

Let $\mathfrak G_{\alpha,e}$ be the stack of pairs $(L,U)$ with
$\dim L=\alpha$ and $U\subset L$ a subrepresentation of dimension $e$.
The assignment $(L,L^\bullet)\mapsto(L,L^n)$ defines an isomorphism
of stacks over $\M_\alpha$:
\[
 \begin{tikzcd}
 \mathfrak F_{\ii,(\alpha-e,e)}
   \arrow[r,"\sim"] \arrow[d,"\pi_{(\alpha-e,e)}"']
   & \mathfrak G_{\alpha,e} \arrow[d,"\rho_{\alpha,e}"] \\
 \M_\alpha \arrow[r,"\mathrm{id}"'] & \M_\alpha,
 \end{tikzcd}
\]
where $\rho_{\alpha,e}(L,U)=L$. Hence
\begin{equation}\label{cc:relative-gr-test}
 \mathcal K_{(\alpha-e,e)}
 =\pi_{(\alpha-e,e)!}\Q\simeq\rho_{\alpha,e!}\Q.
\end{equation}
Since $R(\alpha-e,e)=\eu{e}{\alpha-e}$, Proposition~\ref{char:point-fibre}
and \eqref{geo:character} give
\begin{equation}\label{cc:coefficient-complex}
 \Y_\ii(A_L)_{(\alpha-e,e)}
 \simeq\RGC(\Gr_e(L),\Q)\bal{\eu{e}{\alpha-e}}.
\end{equation}

For the doubled word $\ii=(\iiO,\iiO)$, identify the exponent lattice
$\mathbb Z^{2n}$ of Subsection~\ref{sec:quantum-tori} with
$\mathbb Z^{Q_0}\oplus\mathbb Z^{Q_0}$ in the order of the two copies
of $\iiO$. On $\mathbb Z^{Q_0}=\bigoplus_{i\in Q_0}\mathbb Z\eps_i$,
let $E=(\eu{\eps_i}{\eps_j})$ be the matrix of the Euler form. Set
\begin{equation}\label{cc:matrices}
 B=E^T-E,\qquad
 \widetilde B=\begin{pmatrix}B\\I_n\end{pmatrix},\qquad
 J_\ii=\begin{pmatrix}-E&-E^T\\I_n&0\end{pmatrix}.
\end{equation}
In the target lattice, the first $n$ coordinates are mutable and
the last $n$ are frozen.

Recall the skew-symmetric form $L_\ii$ defining multiplication in $\TT_\ii$.
Define its transported form $\Lambda$ by
\begin{equation}\label{cc:skew-form}
 \Lambda(J_\ii\mathbf a,J_\ii\mathbf b)
 =L_\ii(\mathbf a,\mathbf b),
 \qquad \mathbf a,\mathbf b\in\mathbb Z^{2n}.
\end{equation}
In the standard bases,
\[
 L_\ii=\begin{pmatrix}0&-E^T\\E&0\end{pmatrix},\qquad
 \Lambda=\begin{pmatrix}0&-I_n\\I_n&-B\end{pmatrix},\qquad
 J_\ii^T\Lambda J_\ii=L_\ii.
\]
We use the same notation for a bilinear form and its matrix in the
standard basis, so $\Lambda(g,h)=g^T\Lambda h$.
The change of exponents induces an isomorphism of quantum tori
\begin{equation}\label{cc:torus-isomorphism}
 \widehat J_\ii:\TT_\ii\xrightarrow{\sim}\TT_X,
 \qquad t^{\mathbf a}\longmapsto X^{J_\ii\mathbf a},\qquad
 X^gX^h=w^{\Lambda(g,h)}X^{g+h}.
\end{equation}
The pair $(\widetilde B,\Lambda)$ has principal coefficients and satisfies
the compatibility condition
$\widetilde B^T\Lambda=(I_n\ \ 0)$.
Since $B=E^T-E$, the character exponents are
\[
 J_\ii(\alpha-e,e)
 =(-E(\alpha-e)-E^Te,\alpha-e)
 =(-E\alpha-Be,\alpha-e).
\]

\begin{corollary}[Quantum cluster character]
\label{cc:character}
For every representation $L$ of dimension $\alpha$,
\begin{equation}\label{cc:formula}
 \widehat J_\ii(X_{L,\ii})
 =\sum_{e\le\alpha}w^{-\eu{e}{\alpha-e}}
 |\Gr_e(L)(k)|\,X^{(-E\alpha-Be,\alpha-e)}.
\end{equation}
This is the quantum cluster character for the compatible pair
$(\widetilde B,\Lambda)$, with $w=-q^{1/2}$.
\end{corollary}
\begin{proof}
Take traces in \eqref{cc:coefficient-complex} and apply
$\widehat J_\ii$ to \eqref{char:Z-trace}. The paired shift contributes
$(-1)^{\eu{e}{\alpha-e}}q^{-\eu{e}{\alpha-e}/2}
=w^{-\eu{e}{\alpha-e}}$.
\end{proof}

The realization used in the multiplication formula below is
\begin{equation}\label{cc:Hall-character-normalization}
 \Phi=\widehat J_\ii\circ\Psi_\ii,\qquad
 \Phi([A_L])=\widehat J_\ii(X_{L,\ii}).
\end{equation}

Write $X_L=\widehat J_\ii(X_{L,\ii})$.
For rigid $L$, the quiver
Grassmannians have counting polynomials $P_{L,e}(t)\in\mathbb Z[t]$ with
\[
 |\Gr_e(L)(\mathbb F_{q^r})|=P_{L,e}(q^r),\qquad r\ge1.
\]
Replacing
$|\Gr_e(L)(k)|$ in \eqref{cc:formula} by $P_{L,e}(w^2)$ and regarding
$w$ as an indeterminate gives the generic character, also denoted by $X_L$.

\begin{corollary}[Rigid representations]\label{cc:rigid-specialization}
The assignment $L\mapsto X_L$ gives a bijection from the isomorphism
classes of indecomposable rigid representations to the non-initial quantum
cluster variables for $(\widetilde B,\Lambda)$.
\end{corollary}

\begin{proof}
{With the right-module convention of \cite[\S3]{Qin2012}, taking annihilators
identifies $\Gr_e(L)$ with $\Gr_{\alpha-e}(D_kL)$. Since
$(-E\alpha-Be,\alpha-e)=(-E^T\alpha,0)+\widetilde B(\alpha-e)$,
the assertion follows from \cite[Theorems~3.2.5--3.2.6]{Qin2012} applied to the generic form of \eqref{cc:formula}.}
\end{proof}

\begin{remark}[Classical character]\label{cc:classical-specialization}
For a complex representation $L$ of dimension $\alpha$, the classical character with principal coefficients is
\[
 \sum_{e\le\alpha}\chi_c(\Gr_e(L))X^{(-E\alpha-Be,\alpha-e)}.
\]
For rigid $L$, this is the specialization of the generic character \eqref{cc:formula} at $w=1$.
Setting the frozen variables to $1$ gives the Caldero--Chapoton
formula \cite{CC2006}.
\end{remark}

\section{Triangles for extension spaces of arbitrary dimension}
\label{cc:sheaf-triangle}

We construct an extension triangle for two representations $M,N$.
Put $V=\Ext^1_Q(M,N)$ and assume $r=\dim V>0$.
Starting with a triangle on $[V/\mathbb G_m]$, we push it forward along
the map sending an extension to its middle term. The middle term of
the resulting triangle has a finite filtration with layers given by
shifts and Tate twists of
the Hall convolution $A_M\star A_N$. Pairing with Lusztig sheaves
and taking Frobenius traces produces the factor $1+q+\cdots+q^{r-1}$.
Write $\mathbb L=\Q[-2](-1)$ and
$\mathbb L^s=\Q[-2s](-s)$ for $s\in\mathbb Z$, so that
\[
 \RGC(\mathbb A^s,\Q)\simeq\mathbb L^s
 \qquad(s\geq0).
\]
Indeed, $H_c^2(\mathbb A^1_{\bar k},\Q)=\Q(-1)$ is the only nonzero
cohomology group. The shift $[-2]$ places it in degree $2$, and the
twist $(-1)$ records the geometric Frobenius eigenvalue $q$.

\subsection{A triangle on \texorpdfstring{$[V/\mathbb G_m]$}{[V/Gm]}}

Let $V$ be an $r$-dimensional vector space, $r\geq1$, with the scalar
action of $\mathbb G_m$. Put $\mathfrak U_V=[V/\mathbb G_m]$, with
the open and closed immersions
\[
 j:\mathbb P(V)\simeq[(V\setminus\{0\})/\mathbb G_m]
   \hookrightarrow\mathfrak U_V,\qquad
 i:B\mathbb G_m=[\{0\}/\mathbb G_m]
   \hookrightarrow\mathfrak U_V.
\]
Consider the rank-$r$ vector bundle
\[
 p:\mathscr E=[V\times V/\mathbb G_m]\longrightarrow\mathfrak U_V,
 \qquad [(x,v)]\longmapsto[x],
\]
where $\mathbb G_m$ acts on $V\times V$ by $t\cdot(x,v)=(tx,tv)$.
Its zero section $z$ and diagonal section $s$ are
\[
 z,s:\mathfrak U_V\longrightarrow\mathscr E,
 \qquad z([x])=[(x,0)],\qquad s([x])=[(x,x)].
\]
The complement $\mathscr E^\times$ of the zero section $z$ has the projection
\[
 b:\mathscr E^\times=[V\times(V\setminus\{0\})/\mathbb G_m]
 \longrightarrow\mathfrak U_V,
 \qquad [(x,v)]\longmapsto[x].
\]
The restriction of $s$ to $\mathbb P(V)$ gives
\[
 s^\times=s|_{\mathbb P(V)}:\mathbb P(V)\longrightarrow\mathscr E^\times,
 \qquad [x]\longmapsto[(x,x)],\qquad b\circ s^\times=j.
\]
Set $\mathscr T_V=b_!\Q_{\mathscr E^\times}$.
The diagonal section $s$ vanishes precisely on $i(B\mathbb G_m)$.
We use this to relate the localization triangle for the zero section of
$\mathscr E$ to the two localization triangles for $i$ and $j$.

\begin{proposition}[Triangle and filtration for $\mathscr T_V$]\label{cc:universal-general}
There is a distinguished triangle
\begin{equation}\label{cc:universal-general-triangle}
 j_!\Q_{\mathbb P(V)}\otimes\mathbb L^r
 \longrightarrow\mathscr T_V
 \longrightarrow Rj_*\Q_{\mathbb P(V)}[-1]
 \longrightarrow.
\end{equation}
Let $a:V\to\mathfrak U_V$, $x\mapsto[x]$, be the quotient atlas.
A complete flag in $V$ gives a finite filtration of
$\mathscr T_V$ with layers
\begin{equation}\label{cc:auxiliary-layers}
 a_!\Q_V\otimes\mathbb L^s,\qquad 0\leq s<r.
\end{equation}
\end{proposition}

\begin{proof}
For the complementary immersions $i,j$, the standard localization
triangles of a complex $F$ on $\mathfrak U_V$ are
\[
 j_!j^*F\longrightarrow F\longrightarrow i_*i^*F\longrightarrow,
 \qquad
 i_*i^!F\longrightarrow F\longrightarrow Rj_*j^*F\longrightarrow.
\]
For $F=\Q_{\mathfrak U_V}\otimes\mathbb L^r$, we have
\[
 j^*\Q_{\mathfrak U_V}\simeq\Q_{\mathbb P(V)},
 \qquad i^*\Q_{\mathfrak U_V}\simeq\Q_{B\mathbb G_m}.
\]
The first triangle becomes
\[
 j_!\Q_{\mathbb P(V)}\otimes\mathbb L^r
 \longrightarrow\Q_{\mathfrak U_V}\otimes\mathbb L^r
 \xrightarrow{\ u\ }i_*\Q_{B\mathbb G_m}\otimes\mathbb L^r
 \longrightarrow.
\]
Here $u=\eta\otimes\mathrm{id}$
and $\eta:\mathrm{id}\to i_*i^*$ is the unit.
In the second triangle, take $F=\Q_{\mathfrak U_V}$.
Purity for the zero section $i$ of the rank-$r$ bundle $\mathfrak U_V\to B\mathbb G_m$ gives
\[
 i^!\Q_{\mathfrak U_V}\simeq\Q_{B\mathbb G_m}[-2r](-r)
   =\Q_{B\mathbb G_m}\otimes\mathbb L^r.
\]
With $j^*\Q_{\mathfrak U_V}=\Q_{\mathbb P(V)}$, we obtain
\[
 i_*\Q_{B\mathbb G_m}\otimes\mathbb L^r\xrightarrow{\ v\ }\Q_{\mathfrak U_V}
 \longrightarrow Rj_*\Q_{\mathbb P(V)}
 \xrightarrow{\ \partial\ }.
\]
Consider the composite
\[
 \Q_{\mathfrak U_V}\otimes\mathbb L^r
 \xrightarrow{\ u\ } i_*\Q_{B\mathbb G_m}\otimes\mathbb L^r
 \xrightarrow{\ v\ }\Q_{\mathfrak U_V}.
\]
For the zero section $z$ and its open complement
$\kappa:\mathscr E^\times\hookrightarrow\mathscr E$, the standard
localization triangle of a complex $F$ on $\mathscr E$ is
\[
 \kappa_!\kappa^*F\longrightarrow F\longrightarrow z_*z^*F
 \longrightarrow.
\]
Taking $F=\Q_{\mathscr E}$ and using
$\kappa^*\Q_{\mathscr E}=\Q_{\mathscr E^\times}$ and
$z^*\Q_{\mathscr E}=\Q_{\mathfrak U_V}$ gives
\[
 \kappa_!\Q_{\mathscr E^\times}\longrightarrow\Q_{\mathscr E}\longrightarrow z_*\Q_{\mathfrak U_V}
 \longrightarrow.
\]
We apply $p_!$ to this triangle and compute its three terms.
\[
 p_!\kappa_!\Q_{\mathscr E^\times}\simeq b_!\Q_{\mathscr E^\times}=\mathscr T_V,\quad
 p_!\Q_{\mathscr E}\simeq\Q_{\mathfrak U_V}\otimes\mathbb L^r,\quad
 p_!z_*\Q_{\mathfrak U_V}\simeq(pz)_!\Q_{\mathfrak U_V}=\Q_{\mathfrak U_V}.
\]
Here $p$ is a rank-$r$ vector bundle, while $z_*=z_!$ and
$pz=\mathrm{id}_{\mathfrak U_V}$. We obtain
\[
 \mathscr T_V\longrightarrow\Q_{\mathfrak U_V}\otimes\mathbb L^r
 \xrightarrow{\ \varepsilon\ }\Q_{\mathfrak U_V}
 \longrightarrow.
\]
We show that $\varepsilon=v\circ u$. Write
$\rho:\Q_{\mathscr E}\to z_*\Q_{\mathfrak U_V}$ for the unit above,
so that $\varepsilon=p_!\rho$.
Since $p$ is smooth of relative dimension $r$ and $pz=\mathrm{id}_{\mathfrak U_V}$,
\[
 p^!\Q_{\mathfrak U_V}\simeq\Q_{\mathscr E}[2r](r),
 \qquad z^!p^!\Q_{\mathfrak U_V}\simeq\Q_{\mathfrak U_V}.
\]
The counit $z_*z^!p^!\Q_{\mathfrak U_V}\to p^!\Q_{\mathfrak U_V}$
therefore defines the Gysin map $G_z$, and we set
\[
 \gamma:\Q_{\mathscr E}
 \xrightarrow{\ \rho\ }z_*\Q_{\mathfrak U_V}
 \xrightarrow{\ G_z\ }\Q_{\mathscr E}[2r](r).
\]
If $\theta:p_!p^!\Q_{\mathfrak U_V}\to\Q_{\mathfrak U_V}$ is the
counit, compatibility of counits with composition gives
$\theta\circ p_!G_z=\mathrm{id}_{\Q_{\mathfrak U_V}}$.
Thus $\gamma$ corresponds to $\varepsilon$ under $p_!\dashv p^!$.
Using the identification
$p_!\Q_{\mathscr E}\simeq\Q_{\mathfrak U_V}[-2r](-r)$ induced by
the trace counit, the adjunction identities give
\[
 \gamma=p^*\bigl(\varepsilon[2r](r)\bigr),
 \qquad s^*\gamma=\varepsilon[2r](r),
\]
where the second equality uses $ps=\mathrm{id}_{\mathfrak U_V}$.

Let $p_2:\mathscr E\to\mathfrak U_V$ be the second projection.
Both squares in the diagram
\[
 {\begin{tikzcd}[ampersand replacement=\&,column sep=large]
 B\mathbb G_m \arrow[r,"i"] \arrow[d,"i"']
   \& \mathfrak U_V \arrow[d,"z"] \arrow[r] \& B\mathbb G_m \arrow[d,"i"]\\
 \mathfrak U_V \arrow[r,"s"'] \& \mathscr E \arrow[r,"p_2"'] \& \mathfrak U_V
 \end{tikzcd}}
\]
are Cartesian. Proper base change in the left square gives
\[
 s^*z_*\Q_{\mathfrak U_V}\simeq i_*\Q_{B\mathbb G_m},
 \qquad s^*\rho=u[2r](r).
\]
{Smooth base change in the right square gives $G_z=p_2^*(v[2r](r))$.
Since $p_2s=\mathrm{id}$, we obtain $s^*G_z=v[2r](r)$.}
Consequently,
\[
 \varepsilon[2r](r)=s^*\gamma
   =s^*G_z\circ s^*\rho=(v\circ u)[2r](r),
\]
and hence
\[
 \varepsilon=v\circ u.
\]
Apply the octahedral axiom to $u$ and $v$. In the resulting diagram,
the two full rows and the two full columns are distinguished triangles.
\[
{\begin{tikzcd}[ampersand replacement=\&,column sep=1.9em,row sep=1.5em]
 j_!\Q_{\mathbb P(V)}\otimes\mathbb L^r \arrow[r,equal]\arrow[d]
   \& j_!\Q_{\mathbb P(V)}\otimes\mathbb L^r \arrow[d] \& \& \\
 \mathscr T_V \arrow[r]\arrow[d]
   \& \Q_{\mathfrak U_V}\otimes\mathbb L^r \arrow[r,"\varepsilon"]\arrow[d,"u"]
   \& \Q_{\mathfrak U_V} \arrow[r]\arrow[d,equal] \& \mathscr T_V[1]\arrow[d]\\
 Rj_*\Q_{\mathbb P(V)}[-1]\arrow[r,"{-\partial[-1]}"]\arrow[d]
   \& i_*\Q_{B\mathbb G_m}\otimes\mathbb L^r \arrow[r,"v"']\arrow[d]
   \& \Q_{\mathfrak U_V} \arrow[r] \& Rj_*\Q_{\mathbb P(V)}\\
 (j_!\Q_{\mathbb P(V)}\otimes\mathbb L^r)[1]\arrow[r,equal]
   \& (j_!\Q_{\mathbb P(V)}\otimes\mathbb L^r)[1] \& \&
\end{tikzcd}}
\]
The lower full row is the backward rotation of the triangle for $v$.
The first column gives \eqref{cc:universal-general-triangle}.

This construction of the triangle uses no choice of flag.
To construct the filtration, use the projection
\[
 \pi:\mathscr E^\times\longrightarrow\mathbb P(V),
 \qquad [(x,v)]\longmapsto[v].
\]
Choose a complete flag $0=V_0\subset V_1\subset\cdots\subset V_r=V$,
with $\dim V_d=d$. The flag gives a stratification of $\mathbb P(V)$ with strata
\[
 C_s=\mathbb P(V_{s+1})\setminus\mathbb P(V_s)\simeq\mathbb A^s,
 \qquad 0\leq s<r,
\]
where $\mathbb P(V_0)=\varnothing$.
To describe the inverse images, choose a basis $e_1,\ldots,e_r$
adapted to the flag. Each line in $C_s$ has a unique representative
\[
 v=e_{s+1}+\sum_{d=1}^s u_de_d.
\]
Using the simultaneous scalar action to put the second vector in this
form identifies $\pi^{-1}(C_s)$ with $V\times\mathbb A^s$.
Under this identification, $b|_{\pi^{-1}(C_s)}=a\circ\operatorname{pr}_V$.
Consequently,
\[
 (b|_{\pi^{-1}(C_s)})_!\Q_{\pi^{-1}(C_s)}
 \simeq a_!\Q_V\otimes\RGC(\mathbb A^s,\Q)
 \simeq a_!\Q_V\otimes\mathbb L^s.
\]
Put $X_s=\pi^{-1}(\mathbb P(V_{s+1}))$ for $0\leq s<r$
and $X_{-1}=\varnothing$. Their open complements give
\[
 U_i=\mathscr E^\times\setminus X_{r-1-i},\qquad
 \varnothing=U_0\subset U_1\subset\cdots\subset U_r
 =\mathscr E^\times,
\]
with $U_i\setminus U_{i-1}=\pi^{-1}(C_{r-i})$ for $1\leq i\leq r$.
Set $T_i=(b|_{U_i})_!\Q_{U_i}$, so $T_0=0$ and $T_r=\mathscr T_V$.
Apply the first standard localization triangle to $F=\Q_{U_i}$,
with open substack $U_{i-1}$ and closed complement
$\pi^{-1}(C_{r-i})$, and push forward along $b|_{U_i}$.
The calculation of the stratum contributions above then gives
the required filtration:
\[
 T_{i-1}\longrightarrow T_i\longrightarrow
 a_!\Q_V\otimes\mathbb L^{r-i}\longrightarrow T_{i-1}[1],
 \qquad 1\leq i\leq r.\qedhere
\]
\end{proof}

\begin{remark}\leavevmode
\begin{enumerate}[label=\textup{(\roman*)}]
\item For $r=1$, $\mathscr E^\times\simeq V$ and $b=a$, so
$\mathscr T_V\simeq a_!\Q_V$.
Taking $V=\mathbb A^1$, we can also obtain \eqref{cc:universal-general-triangle}
by applying $a_!$ to the localization triangle at $1$
and using \cite[Lemme~3.2]{Laumon2003}.
For indecomposable rigid $L$ of dimension $\alpha$, let
$j_\alpha:B\mathbb G_m\simeq\mathcal U_\alpha\hookrightarrow\M_\alpha$
denote its open orbit. Pulling back $a$ along $i$ gives the atlas $\pt\to B\mathbb G_m$,
whose composite with $j_\alpha$ is $x_L$. Applying $j_{\alpha!}i^*$ to the zero-section
localization triangle for $\mathscr T_V$ gives
\[
 {A_L\longrightarrow j_{\alpha!}\Q[-2](-1)
 \longrightarrow j_{\alpha!}\Q\longrightarrow A_L[1].}
\]
{Taking the cone of the middle arrow and shifting by $[-1]$ constructs $A_L$ from the rigid open orbit. By \eqref{cc:Hall-character-normalization}, $\Phi([A_L])$ realizes the variable $X_L$ of Corollary~\ref{cc:rigid-specialization} at $w=-q^{1/2}$.}
\item For $r>1$, the filtration \eqref{cc:auxiliary-layers} need not split.
For a geometric point $\bar0$ {of $[\{0\}/\mathbb G_m]$},
\[
 (\mathscr T_V)_{\bar0}\simeq\RGC(\mathbb A^r\setminus\{0\},\Q).
\]
This stalk has nonzero cohomology only in degrees $1$ and $2r$, whereas
\[
 H^d\!\left(\bigoplus_{s=0}^{r-1}
 (a_!\Q_V\otimes\mathbb L^s)_{\bar0}\right)\ne0
 \qquad(1\leq d\leq 2r).
\]
\end{enumerate}
\end{remark}

\subsection{The extension triangle and Hall convolution}
\label{cc:sec:extension-hall}

Fix representations $M,N$, and write
\[
 \begin{gathered}
 m=\dim M,\qquad n=\dim N,\qquad h_{M,N}=\dim\Hom_Q(M,N),\\
 V=\Ext^1_Q(M,N),\qquad r=\dim V>0.
 \end{gathered}
\]
Thus $\eu{m}{n}=h_{M,N}-r$. We construct a family of middle terms $E$
of extensions $0\to N\to E\to M\to0$. Put
\[
 C^0(M,N)=\bigoplus_{i\in Q_0}\Hom_k(M_i,N_i),\qquad
 C^1(M,N)=\bigoplus_{h\in Q_1}\Hom_k(M_{s(h)},N_{t(h)}).
\]
For $z=(z_h)\in C^1(M,N)$, the matrices
\[
 E(z)_i=N_i\oplus M_i,\qquad
 E(z)_h=\begin{pmatrix}N_h&z_h\\0&M_h\end{pmatrix}
\]
define an extension $0\to N\to E(z)\to M\to0$.
Thus $C^1(M,N)$ parameterizes extensions with chosen splittings.
Changing these splittings by
$f=(f_i)\in C^0(M,N)$ changes $z_h$ by
\[
 z_h\longmapsto z_h+N_hf_{s(h)}-f_{t(h)}M_h.
\]
This defines the linear map
\[
 d_{M,N}:C^0(M,N)\longrightarrow C^1(M,N),\qquad
 d_{M,N}(f)_h=N_hf_{s(h)}-f_{t(h)}M_h.
\]
Hence
\[
 V=\Ext^1_Q(M,N)\simeq C^1(M,N)/\operatorname{im}d_{M,N}.
\]
Choose a basis $v_1,\ldots,v_r$ of $V$ and representatives
$z^{(1)},\ldots,z^{(r)}\in C^1(M,N)$. They define a linear section
of the quotient map by
\[
 \sigma:V\longrightarrow C^1(M,N),\qquad
 \sigma\!\left(\sum_{j=1}^r x_jv_j\right)
 =\sum_{j=1}^r x_jz^{(j)}.
\]
Set $E_x=E(\sigma(x))$. Since the entries of the arrow matrices are
regular functions of $x$, they define a family of representations on
the trivial vertex bundles $(N_i\oplus M_i)\otimes_k\mathcal O_V$,
where $\mathcal O_V$ is the structure sheaf of $V$.
For $c\in\mathbb G_m$, set
\[
 g_{c,i}=\begin{pmatrix}c\,\mathrm{id}_{N_i}&0\\0&\mathrm{id}_{M_i}\end{pmatrix}.
\]
Since $g_{c,t(h)}(E_x)_h=(E_{cx})_h g_{c,s(h)}$, these maps give
compatible isomorphisms $E_x\simeq E_{cx}$, hence a morphism
\[
 \phi:\mathfrak U_V\longrightarrow\M_{m+n}.
\]
If $\sigma'$ is another linear section, choose a linear map
$f:V\to C^0(M,N)$ with $\sigma'-\sigma=d_{M,N}\circ f$.
The matrices
\[
 H_i(x)=\begin{pmatrix}\mathrm{id}_{N_i}&-f_i(x)\\0&\mathrm{id}_{M_i}\end{pmatrix}
\]
define an isomorphism of families
$H(x):E_x=E(\sigma(x))\xrightarrow{\sim}E'_x=E(\sigma'(x))$.
Since $f(cx)=cf(x)$, we have
\[
 H_i(cx)g_{c,i}=g_{c,i}H_i(x)\qquad(c\in\mathbb G_m).
\]
Set
\begin{equation}\label{cc:general-objects}
 \begin{gathered}
 \mathcal E_{M,N}=(\phi j)_!\Q,\qquad
 \mathcal B_{M,N}=\phi_!Rj_*\Q[-1],\\
 \mathcal H_{M,N}=\phi_!\mathscr T_V[2r](r).
 \end{gathered}
\end{equation}
Here $\phi j:\mathbb P(V)\to\M_{m+n}$ sends $[x]$ to the middle
term $E_x$ of the corresponding nonsplit extension.
By the projection formula, pairing $\mathcal E_{M,N}$ with $\mathcal K$
gives $\RGC(\mathbb P(V),(\phi j)^*\mathcal K)$.

\begin{theorem}[Extension triangle and convolution filtration]
\label{cc:sheaf-exchange}
For arbitrary $M,N$ with $r>0$, there is a distinguished triangle
\begin{equation}\label{cc:sheaf-exchange-triangle}
 \mathcal E_{M,N}
 \longrightarrow\mathcal H_{M,N}
 \longrightarrow\mathcal B_{M,N}[2r](r)
 \longrightarrow.
\end{equation}
The middle term has a finite filtration with layers
\begin{equation}\label{cc:middle-filtration}
 (A_M\star A_N)\otimes\mathbb L^s,\qquad 0\leq s<r.
\end{equation}
\end{theorem}

\begin{proof}
Push \eqref{cc:universal-general-triangle} along $\phi$ and apply
the paired shift $[2r](r)$. Since
$\mathbb L^r[2r](r)=\Q$, this proves the triangle.

Let $\mathfrak E(M,N)$ be the stack of extensions
$0\to N\to E\to M\to0$ with fixed endpoints $M,N$.
Recall that $q$ in \eqref{eq:hall-correspondence} sends an extension
to its quotient and subrepresentation. With $\xi=(x_M,x_N)$,
we have the cartesian square
\[
 \begin{tikzcd}
 \mathfrak E(M,N) \arrow[r,"\iota"] \arrow[d,"t"']
   & \mathfrak E_{m,n} \arrow[d,"q"] \arrow[r,"p"]
   & \M_{m+n} \\
 \pt \arrow[r,"\xi"'] & \M_m\times\M_n.
 \end{tikzcd}
\]
The composite $\mu=p\iota$ sends an extension to its middle term $E$.
Since $A_M\boxt A_N\simeq\xi_!\Q_{\pt}$, base change in this square gives
\[
 A_M\diamond A_N
 =p_!q^*(A_M\boxt A_N)
 \simeq p_!\iota_!t^*\Q_{\pt}
 \simeq\mu_!\Q_{\mathfrak E(M,N)}.
\]
Put $H=\Hom_Q(M,N)$, viewed as an additive group. The standard
description of the extension stack (cf.~\cite[\S4.5]{Bridgeland2017}
over $\mathbb C$) gives
\[
 \mathfrak E(M,N)\simeq[C^1(M,N)/C^0(M,N)]
 \simeq[V/H]\simeq V\times BH,
\]
where $C^0(M,N)$ acts by $f\cdot z=z+d_{M,N}(f)$ and $H$ acts
trivially on $V$. The same construction works over $k$.
Under the equivalence determined by $\sigma$, the chosen family defines
\[
 \rho:V\longrightarrow\mathfrak E(M,N),\qquad
 \rho=\mathrm{id}_V\times\bigl(\pt\to B\Hom_Q(M,N)\bigr),
 \qquad\mu\rho=\phi a.
\]
The morphism $\rho$ is smooth with geometric fibres
$\Hom_Q(M,N)\simeq\mathbb A^{h_{M,N}}$. Smooth purity identifies the counit
of $\rho_!\dashv\rho^!$ with the trace map
\[
 \rho_!\rho^!\Q_{\mathfrak E(M,N)}
 \simeq\rho_!\Q_V[2h_{M,N}](h_{M,N})
 \xrightarrow{\operatorname{Tr}_\rho}\Q_{\mathfrak E(M,N)}.
\]
At a geometric point $\bar e$ of $\mathfrak E(M,N)$, base change
identifies its stalk with the standard isomorphism
\[
 \bigl(\rho_!\Q_V[2h_{M,N}](h_{M,N})\bigr)_{\bar e}
 \simeq R\Gamma_c(\mathbb A^{h_{M,N}},\Q)[2h_{M,N}](h_{M,N})
 \xrightarrow[\sim]{(\operatorname{Tr}_\rho)_{\bar e}}\Q.
\]
Thus $\operatorname{Tr}_\rho$ is an isomorphism, and
\[
 \rho_!\Q_V\simeq\Q_{\mathfrak E(M,N)}\otimes\mathbb L^{h_{M,N}}.
\]
Using $\mu\rho=\phi a$ and the projection formula, we obtain
\[
 \begin{aligned}
 (\phi a)_!\Q_V
 &\simeq\mu_!\rho_!\Q_V
 \simeq\mu_!\bigl(\Q_{\mathfrak E(M,N)}\otimes\mathbb L^{h_{M,N}}\bigr)
 \\
 &\simeq\bigl(\mu_!\Q_{\mathfrak E(M,N)}\bigr)\otimes\mathbb L^{h_{M,N}}
 \simeq(A_M\diamond A_N)\otimes\mathbb L^{h_{M,N}}.
 \end{aligned}
\]
Hence each layer of
$\phi_!\mathscr T_V[2r](r)$ is
\[
 (A_M\diamond A_N)\otimes\mathbb L^{h_{M,N}+s}[2r](r)
 =(A_M\diamond A_N)[2(r-h_{M,N})](r-h_{M,N})\otimes\mathbb L^s
 =(A_M\star A_N)\otimes\mathbb L^s,
\]
as required.
\end{proof}

For $r=1$, $\mathcal H_{M,N}\simeq A_M\star A_N$.
If $\Ext^1_Q(M,N)=0$, the fixed-endpoint extension stack is $BH$,
where $H=\Hom_Q(M,N)$. The same trace isomorphism for $\pt\to BH$
gives $A_M\star A_N\simeq A_{M\oplus N}$.

\subsection{Pairing the extension triangle}

For every bounded mixed complex $\mathcal K$ on $\M_{m+n}$, exactness gives
\begin{equation}\label{cc:paired-exchange-triangle}
 \mathsf P(\mathcal E_{M,N},\mathcal K)
 \longrightarrow\mathsf P(\mathcal H_{M,N},\mathcal K)
 \longrightarrow\mathsf P(\mathcal B_{M,N},\mathcal K)[2r](r)
 \longrightarrow.
\end{equation}
The filtration of $\mathcal H_{M,N}$ gives the middle term a finite
filtration with layers
\[
 \mathsf P(A_M\star A_N,\mathcal K)\otimes\mathbb L^s,
 \qquad 0\leq s<r.
\]
For a Lusztig sheaf $\mathcal K=\mathcal K_{\mathbf c}$ with
$\wt{\mathbf c}=m+n$, Lemma~\ref{geo:pairing-transposition} gives
\[
 \mathsf P(A_M\star A_N,\mathcal K_{\mathbf c})
 \simeq(x_M,x_N)^*\mathcal R_{m,n}(\mathcal K_{\mathbf c})[2(r-h_{M,N})](r-h_{M,N}).
\]
By Proposition~\ref{geo:restriction-filtration}, each $s$-layer therefore
has a further filtration with layers
\[
 \bigl(x_M^*\mathcal K_{\mathbf a}\otimes x_N^*\mathcal K_{\mathbf b}\bigr)
 [2\ef{\mathbf a}{\mathbf b}+2(r-h_{M,N})](\ef{\mathbf a}{\mathbf b}+r-h_{M,N})\otimes\mathbb L^s,
\]
where $\mathbf a+\mathbf b=\mathbf c$, $\wt{\mathbf a}=m$, and
$\wt{\mathbf b}=n$. After normalization
\eqref{geo:character}, these are the layers of
Theorem~\ref{geo:main-theorem}, tensored with $\mathbb L^s$.

For any $\mathcal K$ as above, trace additivity gives
\[
 \Tr\bigl(\Fr,\mathsf P(\mathcal H_{M,N},\mathcal K)\bigr)
 =[r]_q\,\Tr\bigl(\Fr,\mathsf P(A_M\star A_N,\mathcal K)\bigr),
 \qquad [r]_q=\sum_{s=0}^{r-1}q^s.
\]
This identity requires no splitting of the filtration.
In Section~\ref{sec:boundary-multiplication}, we identify
the boundary coefficients and derive the general multiplication formula.
\section{Boundary coefficients and general multiplication}
\label{sec:boundary-multiplication}

We compute the coefficients of $\mathcal B_{M,N}$ in the extension
triangle of Theorem~\ref{cc:sheaf-exchange} and express them in terms
of quiver Grassmannians. These calculations yield the general
multiplication formula.

\subsection{The boundary coefficient theorem}
\label{bd:sec:boundary-coefficients}
Fix representations $M,N$ and put
\[
 \begin{gathered}
 m=\dim M,\qquad n=\dim N,\\
 V=\Ext^1_Q(M,N),\qquad r=\dim V>0.
 \end{gathered}
\]
Use the family $E_x=E(\sigma(x))$ constructed before
Theorem~\ref{cc:sheaf-exchange}, with
\[
 \phi:\mathfrak U_V=[V/\mathbb G_m]\longrightarrow\M_{m+n},
 \qquad j:\mathbb P(V)\hookrightarrow\mathfrak U_V.
\]
Recall
\begin{equation}\label{bd:eq:boundary}
 \mathcal B_{M,N}=\phi_!Rj_*\Q[-1].
\end{equation}
Here ``boundary'' refers to the third term
$\mathcal B_{M,N}[2r](r)$ of the extension triangle
\eqref{cc:sheaf-exchange-triangle}.
Put $\alpha=m+n$. For a dimension vector
$e=(e_i)_{i\in Q_0}\in\mathbb N^{Q_0}$ with $e\le\alpha$, recall the projection
\[
 \rho_{\alpha,e}:\mathfrak G_{\alpha,e}\longrightarrow\M_\alpha,
 \qquad(E,L)\longmapsto E,
\]
where $L\subset E$ is a subrepresentation of dimension $e$.
Recall that $\mathcal K_{\mathbf a}=\pi_{\mathbf a!}\Q$, where
$\pi_{\mathbf a}:(E,E^\bullet)\mapsto E$ forgets a flag with
$E^{j-1}/E^j\simeq S_{i_j}^{\oplus a_j}$.
For $\ii=(\iiO,\iiO)$, a flag of type $(\alpha-e,e)$ is obtained
by joining the source-adapted flags of $E/L$ and $L$, of dimensions
$\alpha-e$ and $e$, respectively. By \eqref{cc:relative-gr-test},
\[
 \mathcal K_{(\alpha-e,e)}\simeq\rho_{\alpha,e!}\Q.
\]
Its boundary coefficients are the pairings
\eqref{geo:pairing-definition}:
\[
 \mathsf P(\mathcal B_{M,N},\mathcal K_{(\alpha-e,e)})
 =\RGC(\M_\alpha,\mathcal B_{M,N}\otimes\mathcal K_{(\alpha-e,e)}).
\]

For dimension vectors $u,v\in\mathbb N^{Q_0}$ with
$u\le m$, $v\le n$, and $u+v=e$, set
$S_{u,v}=\Gr_u(M)\times\Gr_v(N)$.
For $(U,W)\in S_{u,v}$ and an extension
\[
 0\longrightarrow N\longrightarrow E_x\xrightarrow{p_x}M\longrightarrow0,
\]
its pullback to $U$ has middle term $p_x^{-1}(U)$. Quotienting by
$W\subset N$ then gives
\[
 0\longrightarrow N/W\longrightarrow p_x^{-1}(U)/W
   \longrightarrow U\longrightarrow0.
\]
Let $a_U$ denote pullback along $U\hookrightarrow M$ and $b_W$
pushout along $N\twoheadrightarrow N/W$. Sending $x$ to the class
of this extension defines their composite
\begin{equation}\label{bd:eq:beta}
 \beta_{U,W}=b_W\circ a_U:
 V\xrightarrow{a_U}\Ext^1_Q(U,N)
                   \xrightarrow{b_W}\Ext^1_Q(U,N/W).
\end{equation}
The short exact sequences $0\to U\to M\to M/U\to0$ and
$0\to W\to N\to N/W\to0$ give
\[
 \begin{aligned}
 \Ext^1_Q(M,N)&\xrightarrow{a_U}\Ext^1_Q(U,N)
                  \longrightarrow\Ext^2_Q(M/U,N)=0,\\
 \Ext^1_Q(U,N)&\xrightarrow{b_W}\Ext^1_Q(U,N/W)
                  \longrightarrow\Ext^2_Q(U,W)=0.
 \end{aligned}
\]
Since the category $\operatorname{rep}_kQ$ is hereditary, every representation has
a projective resolution of length at most one, so $\Ext^i_Q(-,-)=0$
for $i\ge2$.
Exactness shows that $a_U$ and $b_W$ are surjective, hence so is
$\beta_{U,W}=b_W\circ a_U$.

The quotient extension splits precisely when a subrepresentation
$L\subset E_x$ has $L\cap N=W$ and $p_x(L)=U$: the splitting
identifies $U$ with $L/W$. Thus
\[
 \beta_{U,W}(x)=0
 \quad\Longleftrightarrow\quad
 \exists\,L\subset E_x:\ L\cap N=W,\quad p_x(L)=U.
\]
The Fourier calculation in Lemma~\ref{bd:lem:fourier} expresses the
boundary contribution from this kernel in terms of the linear
functionals that vanish on it. With $D_k=\Hom_k(-,k)$, we therefore define
\begin{equation}\label{bd:eq:xi}
 \Xi_{u,v}=
 \{(U,W,[\eta])\in S_{u,v}\times\mathbb P(D_kV):
                     \eta|_{\ker\beta_{U,W}}=0\}.
\end{equation}
Dualizing the short exact sequence for $\beta_{U,W}$ gives
\[
 0\longrightarrow D_k\Ext^1_Q(U,N/W)
 \xrightarrow{D_k\beta_{U,W}}D_kV
 \xrightarrow{\eta\mapsto\eta|_{\ker\beta_{U,W}}}
 D_k\ker\beta_{U,W}\longrightarrow0.
\]
Hence
\[
 (\Xi_{u,v})_{(U,W)}
 \simeq\mathbb P\bigl(D_k\Ext^1_Q(U,N/W)\bigr).
\]
We give its closed subscheme structure in
Subsection~\ref{bd:sec:lifting}.

\begin{theorem}[Boundary coefficients]\label{bd:thm:main}
For every $e$, the complex
$\mathsf P(\mathcal B_{M,N},\mathcal K_{(\alpha-e,e)})$ has a finite
filtration with layers
\begin{equation}\label{bd:eq:layers}
 \RGC(\Xi_{u,v},\Q)[-2t_{u,v}](-t_{u,v}),\qquad
 t_{u,v}=r+\eu{u}{n-v},\qquad u+v=e.
\end{equation}
For each nonempty $S_{u,v}$, one has $t_{u,v}\ge0$.
Consequently, its Frobenius trace is
\begin{equation}\label{bd:eq:trace}
 \begin{aligned}
 b_e(q)&:=\Tr\bigl(\Fr,\mathsf P(\mathcal B_{M,N},\mathcal K_{(\alpha-e,e)})\bigr)\\
       &=\sum_{u+v=e}q^{\,r+\eu{u}{n-v}}\#\Xi_{u,v}(k).
 \end{aligned}
\end{equation}
The same formula holds over finite extensions of $k$.
\end{theorem}

With $R_e=\eu{e}{m+n-e}$, the normalized coefficient
$\mathsf P(\mathcal B_{M,N},\mathcal K_{(\alpha-e,e)})\bal{R_e}$
of Section~\ref{sec:characters} therefore has layers
\begin{equation}\label{bd:eq:normalized}
 \RGC(\Xi_{u,v},\Q)\bal{R_e-2t_{u,v}}.
\end{equation}
The trace formula uses additivity for this filtration and does not
require a splitting or polynomial point counts for $\Xi_{u,v}$.

\subsection{Fourier transform and dual zero loci}\label{bd:sec:fourier}

We give the Fourier calculation in detail to show how $\Xi_{u,v}$ enters the multiplication formula.
In Subsection~\ref{bd:sec:lifting}, we apply this calculation to the subrepresentation
strata in the boundary coefficient.

Let $S$ be a finite-type $k$-scheme, let $\pi_{\mathcal V}:\mathcal V\to S$
be a vector bundle of rank $r\ge1$, and let
$\pi_{\mathcal V^\vee}:\mathcal V^\vee\to S$ be its dual bundle.
Write
\[
 i_0:S\to\mathcal V,\qquad i_0^\vee:S\to\mathcal V^\vee
\]
for the zero sections, whose images we denote by $0$, and
\[
 j_0:\mathcal V\setminus0\hookrightarrow\mathcal V,\qquad
 j_0^\vee:\mathcal V^\vee\setminus0\hookrightarrow\mathcal V^\vee
\]
for the complementary open immersions.

Let $p_1,p_2$ be the projections from $\mathcal V\times_S\mathcal V^\vee$
to $\mathcal V$ and $\mathcal V^\vee$, respectively, and let
\[
 \operatorname{ev}:\mathcal V\times_S\mathcal V^\vee\longrightarrow\mathbb A^1_S,
 \qquad (v,\xi)\longmapsto(b,\xi(v)),
\]
where $v\in\mathcal V_b$ and $\xi\in\mathcal V_b^\vee$ lie over the same point
$b\in S$.

Fix a nontrivial additive character
$\psi:(\mathbb F_p,+)\to\overline{\mathbb Q}_\ell^\times$.
Its Artin--Schreier sheaf $\mathcal L_\psi$ is a rank-one local system on
$\mathbb A^1_k$ with Frobenius trace
\[
 \operatorname{Tr}(\operatorname{Fr}_x,(\mathcal L_\psi)_{\bar x})
 =\psi\bigl(\operatorname{Tr}_{k/\mathbb F_p}(x)\bigr),\qquad x\in k.
\]
We also write $\mathcal L_\psi$ for its pullback to $\mathbb A^1_S$.
Use the Fourier transform
\[
 \operatorname{Four}_\psi(K)
 =Rp_{1!}\bigl(p_2^*K\otimes\operatorname{ev}^*\mathcal L_\psi\bigr).
\]

\begin{lemma}\label{bd:lem:fourierkernel}
With this notation, there is an isomorphism
\begin{equation}\label{bd:eq:kernel}
 \operatorname{Four}_\psi(j_{0!}^\vee\Q)\simeq Rj_{0*}\Q[-1].
\end{equation}
It is compatible with the $\mathbb G_m$-action
$c\cdot(v,\xi)=(cv,c^{-1}\xi)$ on $\mathcal V\times_S\mathcal V^\vee$.
\end{lemma}

\begin{proof}
Start with the localization triangle
\[
 j_{0!}^\vee\Q\longrightarrow\Q_{\mathcal V^\vee}
                   \xrightarrow{u}i_{0*}^\vee\Q_S\longrightarrow.
\]
Here $u$ is the unit
$\operatorname{id}\to i_{0*}^\vee(i_0^\vee)^*$.
The following two Fourier transforms are special cases of \cite[Proposition~1.2.3.1]{Laumon1987}. We give the calculations in our convention.
For a geometric point $v\in\mathcal V_b$, write
\[
 \operatorname{ev}_v:\mathcal V_b^\vee\longrightarrow\mathbb A^1,
 \qquad\xi\longmapsto\xi(v).
\]
Then
\[
 \bigl(\operatorname{Four}_\psi(\Q_{\mathcal V^\vee})\bigr)_v
 \simeq \bigl(Rp_{1!}\operatorname{ev}^*\mathcal L_\psi\bigr)_v
 \simeq R\Gamma_c\bigl(\mathcal V_b^\vee,
                         \operatorname{ev}_v^*\mathcal L_\psi\bigr).
\]
If $v=0$, the map $\operatorname{ev}_0$ factors through the point $0\in\mathbb A^1$.
Its pullback is therefore constant with fibre $(\mathcal L_\psi)_0\simeq\Q$, and
\[
 R\Gamma_c(\mathcal V_b^\vee,\operatorname{ev}_0^*\mathcal L_\psi)
 \simeq R\Gamma_c(\mathbb A^r,\Q)\simeq\Q[-2r](-r).
\]
If $v\ne0$, choose $\xi_0\in\mathcal V_b^\vee$ with $\xi_0(v)=1$.
The isomorphism
\[
 \ker(\operatorname{ev}_v)\times\mathbb A^1\xrightarrow{\sim}\mathcal V_b^\vee,
 \qquad(\zeta,a)\longmapsto\zeta+a\xi_0
\]
identifies $\operatorname{ev}_v$ with the second projection, and
$\ker(\operatorname{ev}_v)\simeq\mathbb A^{r-1}$. Thus
$\operatorname{ev}_v^*\mathcal L_\psi\simeq\Q\boxtimes\mathcal L_\psi$, and
the K\"unneth formula gives
\[
 R\Gamma_c(\mathcal V_b^\vee,\operatorname{ev}_v^*\mathcal L_\psi)
 \simeq R\Gamma_c(\mathbb A^{r-1},\Q)\otimes
 R\Gamma_c(\mathbb A^1,\mathcal L_\psi)=0,
\]
since $R\Gamma_c(\mathbb A^1,\mathcal L_\psi)=0$ for nontrivial $\psi$.
Thus $\operatorname{Four}_\psi(\Q_{\mathcal V^\vee})$ is supported on $i_0(S)$.
The map $\widetilde i_0(\xi)=(0_b,\xi)$ gives the Cartesian square
\[
 \begin{tikzcd}
 \mathcal V^\vee \arrow[r,"\widetilde i_0"]
       \arrow[d,"\pi_{\mathcal V^\vee}"']
   & \mathcal V\times_S\mathcal V^\vee \arrow[d,"p_1"] \\
 S \arrow[r,"i_0"'] & \mathcal V .
 \end{tikzcd}
\]
Since $\operatorname{ev}(0_b,\xi)=(b,0)$,
$\widetilde i_0^*\operatorname{ev}^*\mathcal L_\psi\simeq\Q_{\mathcal V^\vee}$.
Thus
\[
 \begin{aligned}
 i_0^*\operatorname{Four}_\psi(\Q_{\mathcal V^\vee})
 &\simeq i_0^*Rp_{1!}\operatorname{ev}^*\mathcal L_\psi\\
 &\simeq R\pi_{\mathcal V^\vee!}
                   \widetilde i_0^*\operatorname{ev}^*\mathcal L_\psi\\
 &\simeq R\pi_{\mathcal V^\vee!}\Q_{\mathcal V^\vee}
 \simeq\Q_S[-2r](-r).
 \end{aligned}
\]
The last identification is induced by the trace isomorphism
\[
 R\pi_{\mathcal V^\vee!}\pi_{\mathcal V^\vee}^!\Q_S
 \simeq R\pi_{\mathcal V^\vee!}\Q_{\mathcal V^\vee}[2r](r)
 \xrightarrow[\sim]{\operatorname{Tr}_{\pi_{\mathcal V^\vee}}}\Q_S.
\]
Consequently,
\[
 \operatorname{Four}_\psi(\Q_{\mathcal V^\vee})
 \simeq i_{0*}\Q_S[-2r](-r).
\]
The pullback $p_2^*i_{0*}^\vee\Q_S$ is supported on
$\mathcal V\times_S i_0^\vee(S)\simeq\mathcal V$, where
\[
 \operatorname{ev}(v,0)=(b,0),\qquad
 (\operatorname{ev}^*\mathcal L_\psi)|_{\mathcal V\times_S i_0^\vee(S)}
 \simeq\Q_{\mathcal V},\qquad
 p_1|_{\mathcal V\times_S i_0^\vee(S)}=\operatorname{id}_{\mathcal V}.
\]
Hence
\[
 \operatorname{Four}_\psi(i_{0*}^\vee\Q_S)
 =Rp_{1!}\bigl(p_2^*i_{0*}^\vee\Q_S\otimes\operatorname{ev}^*\mathcal L_\psi\bigr)
 \simeq R(\operatorname{id}_{\mathcal V})_!\Q_{\mathcal V}
 \simeq\Q_{\mathcal V}.
\]
Applying the Fourier transform to the localization triangle gives
\[
 \operatorname{Four}_\psi(j_{0!}^\vee\Q)
 \longrightarrow i_{0*}\Q_S[-2r](-r)
 \xrightarrow{\operatorname{Four}_\psi(u)}\Q_{\mathcal V}\longrightarrow.
\]
We compare $\operatorname{Four}_\psi(u)$ with the Gysin map for $i_0$.
Applying the preceding calculation of
$\operatorname{Four}_\psi(\Q_{\mathcal V^\vee})$ to the dual bundle gives
\[
 Rp_{2!}\operatorname{ev}^*\mathcal L_\psi
 \simeq i_{0*}^\vee\Q_S[-2r](-r).
\]
Since $\pi_{\mathcal V}p_1=\pi_{\mathcal V^\vee}p_2$, we have
\[
 \begin{aligned}
 R\pi_{\mathcal V!}\operatorname{Four}_\psi(K)
 &\simeq R\pi_{\mathcal V^\vee!}Rp_{2!}
       \bigl(p_2^*K\otimes\operatorname{ev}^*\mathcal L_\psi\bigr)\\
 &\simeq R\pi_{\mathcal V^\vee!}
       \bigl(K\otimes Rp_{2!}\operatorname{ev}^*\mathcal L_\psi\bigr)\\
 &\simeq R\pi_{\mathcal V^\vee!}
       \bigl(K\otimes i_{0*}^\vee\Q_S\bigr)[-2r](-r)\\
 &\simeq R\pi_{\mathcal V^\vee!}i_{0*}^\vee(i_0^\vee)^*K[-2r](-r)\\
 &\simeq(i_0^\vee)^*K[-2r](-r).
 \end{aligned}
\]
Since $(i_0^\vee)^*u=\mathrm{id}_{\Q_S}$, functoriality in $K$ shows that
$R\pi_{\mathcal V!}\operatorname{Four}_\psi(u)$ is an isomorphism.
Now $\pi_{\mathcal V}^!\Q_S\simeq\Q_{\mathcal V}[2r](r)$ and
$\pi_{\mathcal V}i_0=\operatorname{id}_S$ give the adjunction
\[
 \operatorname{Hom}_{\mathcal V}
   (i_{0*}\Q_S,\pi_{\mathcal V}^!\Q_S)
 \simeq\operatorname{End}_S(\Q_S),\qquad
 a\longmapsto\operatorname{Tr}_{\pi_{\mathcal V}}\circ R\pi_{\mathcal V!}a.
\]
Under this adjunction, $\mathrm{Gys}_{i_0}[2r](r)$ corresponds to
$\mathrm{id}_{\Q_S}$ by composition of counits. The shifted Fourier map
corresponds to an automorphism $b$ of $\Q_S$, since the trace map for
the vector bundle $\pi_{\mathcal V}$ is an isomorphism. Thus
\[
 \operatorname{Four}_\psi(u)[2r](r)=\mathrm{Gys}_{i_0}[2r](r)\circ i_{0*}b.
\]
Purity $i_0^!\Q_{\mathcal V}\simeq\Q_S[-2r](-r)$ and the localization
triangle for $i_0,j_0$ give the lower row below.
\[
\begin{tikzcd}[column sep=4em,row sep=2.1em,ampersand replacement=\&]
 \operatorname{Four}_\psi(j_{0!}^\vee\Q)
   \arrow[r]\arrow[d,dashed,"\sim"'] \&
 i_{0*}\Q_S[-2r](-r)
   \arrow[r,"{\operatorname{Four}_\psi(u)}"]
   \arrow[d,"{(i_{0*}b)[-2r](-r)}"] \&
 \Q_{\mathcal V}\arrow[d,"\mathrm{id}"]\\
 Rj_{0*}\Q_{\mathcal V\setminus0}[-1]\arrow[r] \&
 i_{0*}\Q_S[-2r](-r)\arrow[r,"{\mathrm{Gys}_{i_0}}"'] \&
 \Q_{\mathcal V}
\end{tikzcd}
\]
The identity above makes the right square commute. It extends to a morphism
of triangles, and the dashed arrow is an isomorphism since the other two
vertical arrows are isomorphisms.
The same argument over $S\times B\mathbb G_m$, with weights $1$ and $-1$ on
$\mathcal V$ and $\mathcal V^\vee$, gives the equivariant isomorphism.
\end{proof}

Let $\mathcal Z\to S$ be a vector bundle of rank $t$ and let
$f:\mathcal Z\to\mathcal V$ be a morphism of vector bundles.
Use the fibrewise scalar actions of $\mathbb G_m$
and write $\mathbb P(\mathcal V)=[(\mathcal V\setminus0)/\mathbb G_m]$.
Let
\[
 \bar f:[\mathcal Z/\mathbb G_m]\to[\mathcal V/\mathbb G_m],\qquad
 \rho:[\mathcal Z/\mathbb G_m]\to S,\qquad
 j:\mathbb P(\mathcal V)\hookrightarrow[\mathcal V/\mathbb G_m].
\]
Write $f_b=f|_{\mathcal Z_b}:\mathcal Z_b\to\mathcal V_b$ for the fibre map at
$b\in S$. For $\pi:\mathbb P(\mathcal V^\vee)\to S$, set
\[
 \Xi_f=\{(b,[\eta])\in\mathbb P(\mathcal V^\vee):
          \eta\circ f_b=0\},
 \qquad p=\pi|_{\Xi_f}:\Xi_f\to S.
\]
The closed subscheme structure on $\Xi_f$ is defined by the vanishing of
\[
 \mathcal O(-1)\hookrightarrow\pi^*\mathcal V^\vee
     \xrightarrow{\pi^*f^\vee}\pi^*\mathcal Z^\vee,
\]
where $\mathcal O(-1)\subset\pi^*\mathcal V^\vee$ is the tautological line subbundle
whose fibre at $(b,[\eta])$ is the line spanned by $\eta$.

\begin{lemma}\label{bd:lem:fourier}
With this notation, there is an isomorphism
\begin{equation}\label{bd:eq:fourierlemma}
 R\rho_!\bar f^*Rj_*\Q[-1]
 \simeq Rp_*\Q_{\Xi_f}[-2t](-t).
\end{equation}
\end{lemma}

\begin{proof}
The square
\[
 \begin{tikzcd}
 \mathcal Z\times_S\mathcal V^\vee
   \arrow[r,"f\times\operatorname{id}"] \arrow[d,"\widetilde p_1"']
   & \mathcal V\times_S\mathcal V^\vee \arrow[d,"p_1"] \\
 \mathcal Z \arrow[r,"f"'] & \mathcal V
 \end{tikzcd}
\]
is Cartesian. Write $\widetilde p_2=p_2\circ(f\times\operatorname{id})$
for the second projection, put
$\operatorname{pr}_{\mathcal Z}
=\widetilde p_1\big|_{\mathcal Z\times_S(\mathcal V^\vee\setminus0)}$,
and set $\operatorname{ev}_f(z,\xi)=\xi(f(z))$.
By Lemma~\ref{bd:lem:fourierkernel},
\[
 \begin{aligned}
 f^*Rj_{0*}\Q[-1]
 &\simeq f^*\operatorname{Four}_\psi(j_{0!}^\vee\Q)\\
 &=f^*Rp_{1!}\bigl(p_2^*j_{0!}^\vee\Q
              \otimes\operatorname{ev}^*\mathcal L_\psi\bigr)\\
 &\simeq R\widetilde p_{1!}\bigl(\widetilde p_2^*j_{0!}^\vee\Q
              \otimes\operatorname{ev}_f^*\mathcal L_\psi\bigr)\\
 &\simeq R\operatorname{pr}_{\mathcal Z!}
              \operatorname{ev}_f^*\mathcal L_\psi.
 \end{aligned}
\]
Form
\[
 q:\mathfrak Z=[\mathcal Z\times_S(\mathcal V^\vee\setminus0)/\mathbb G_m]
 \longrightarrow\mathbb P(\mathcal V^\vee),\qquad[(z,\xi)]\longmapsto[\xi].
\]
The action $c\cdot(z,\eta)=(cz,c^{-1}\eta)$ gives an isomorphism
\[
 \mathfrak Z\xrightarrow{\sim}\pi^*\mathcal Z\otimes\mathcal O(-1),
 \qquad[(z,\eta)]\longmapsto z\otimes\eta,
\]
since $(cz)\otimes(c^{-1}\eta)=z\otimes\eta$. Locally, choosing a generator
$\xi$ of $\mathcal O(-1)$ gives the inverse
$z\otimes\xi_x\mapsto[(z,\xi_x)]$. Thus $q$ is a rank-$t$ vector bundle.
Since $f$ is linear, $(c^{-1}\xi)(f(cz))=\xi(f(z))$, so
$\operatorname{ev}_f$ descends to
\[
 \ell:\mathfrak Z\longrightarrow\mathbb A^1_S,\qquad
 \ell([(z,\xi)])=(b,\xi(f_b(z))),\qquad
 z\in\mathcal Z_b,\quad\xi\in\mathcal V_b^\vee\setminus0.
\]
This function is linear on the fibres of $q$.

Let $o:\mathbb P(\mathcal V^\vee)\to\mathfrak Z^\vee$
be the zero section of the dual bundle. Since
$(\mathfrak Z^\vee)_{[\xi]}=\bigl(q^{-1}([\xi])\bigr)^\vee$,
$\ell$ defines a section
\[
 s_\ell:\mathbb P(\mathcal V^\vee)\to\mathfrak Z^\vee,
 \qquad s_\ell([\xi])=\ell|_{q^{-1}([\xi])}.
\]
Its zero locus is
\[
 \{(b,[\xi])\in\mathbb P(\mathcal V^\vee):s_\ell([\xi])=0\}
 =\{(b,[\xi])\in\mathbb P(\mathcal V^\vee):\xi\circ f_b=0\}
 =\Xi_f.
\]
Choose a local generator $\xi$ of $\mathcal O(-1)$ on an open subset
$U\subset\mathbb P(\mathcal V^\vee)$. At $x=(b,[\eta])\in U$,
$\xi_x\in\mathcal O(-1)_x\setminus\{0\}$ and $[\xi_x]=[\eta]$.
The generator $\xi$ identifies $\mathcal O_U$ with $\mathcal O(-1)|_U$
by $a\mapsto a\xi$. Tensoring with $\pi^*\mathcal Z|_U$ gives
\[
 \pi^*\mathcal Z|_U\xrightarrow{\sim}\mathfrak Z|_U,
 \qquad(x,z)\longmapsto[(z,\xi_x)].
\]
{Under this trivialization, $s_\ell$ is given by
$(\pi^*f^\vee)(\xi)$. Indeed, for $x=(b,[\eta])$ and $z\in\mathcal Z_b$,}
\[
 s_\ell(x)([(z,\xi_x)])
 =\xi_x(f_b(z))=\bigl(f_b^\vee(\xi_x)\bigr)(z).
\]
Since $\xi$ generates $\mathcal O(-1)$ on $U$, the vanishing of
$(\pi^*f^\vee)(\xi)$ defines $\Xi_f$ on $U$. Hence the zero scheme of $s_\ell$ is $\Xi_f$.
To compute $Rq_!\ell^*\mathcal L_\psi$, first suppose $t>0$.
Apply the preceding Fourier calculation to the rank-$t$ bundle
$q:\mathfrak Z\to\mathbb P(\mathcal V^\vee)$:
\[
 \operatorname{Four}_\psi(\Q_{\mathfrak Z})
 \simeq o_*\Q_{\mathbb P(\mathcal V^\vee)}[-2t](-t).
\]
For $t=0$, the same identity holds since the Fourier transform
on the zero bundle is the identity and $o$ is an isomorphism.
Writing $\iota:\Xi_f\hookrightarrow\mathbb P(\mathcal V^\vee)$, we have the Cartesian square
\[
 \begin{tikzcd}
 \Xi_f \arrow[r,"\iota"] \arrow[d,"\iota"']
   & \mathbb P(\mathcal V^\vee) \arrow[d,"o"] \\
 \mathbb P(\mathcal V^\vee) \arrow[r,"s_\ell"'] & \mathfrak Z^\vee .
 \end{tikzcd}
\]
Let $\operatorname{ev}_{\mathfrak Z}:\mathfrak Z\times_{\mathbb P(\mathcal V^\vee)}\mathfrak Z^\vee
\to\mathbb A^1_S$ be evaluation. The section $s_\ell$ induces
\[
 (\operatorname{id},s_\ell\circ q):\mathfrak Z
 \longrightarrow\mathfrak Z\times_{\mathbb P(\mathcal V^\vee)}\mathfrak Z^\vee,
 \qquad w\longmapsto(w,s_\ell(q(w))).
\]
Since $s_\ell(q(w))(w)=\ell(w)$, we have
\[
 (\operatorname{id},s_\ell\circ q)^*
 \operatorname{ev}_{\mathfrak Z}^*\mathcal L_\psi
 =(\operatorname{ev}_{\mathfrak Z}\circ(\operatorname{id},s_\ell\circ q))^*
 \mathcal L_\psi=\ell^*\mathcal L_\psi.
\]
Let $\operatorname{pr}_2$ be the second projection. The Cartesian square
\[
 \begin{tikzcd}[column sep=large]
 \mathfrak Z \arrow[r,"{(\operatorname{id},s_\ell\circ q)}"] \arrow[d,"q"']
 & \mathfrak Z\times_{\mathbb P(\mathcal V^\vee)}\mathfrak Z^\vee
     \arrow[d,"\operatorname{pr}_2"] \\
 \mathbb P(\mathcal V^\vee) \arrow[r,"s_\ell"'] & \mathfrak Z^\vee
 \end{tikzcd}
\]
gives
\[
 s_\ell^*\operatorname{Four}_\psi(\Q_{\mathfrak Z})
 =s_\ell^*R\operatorname{pr}_{2!}
       \operatorname{ev}_{\mathfrak Z}^*\mathcal L_\psi
 \simeq Rq_!(\operatorname{id},s_\ell\circ q)^*
       \operatorname{ev}_{\mathfrak Z}^*\mathcal L_\psi
 =Rq_!\ell^*\mathcal L_\psi.
\]
Base change in the square defining $\Xi_f$ therefore yields
\[
 Rq_!\ell^*\mathcal L_\psi
 \simeq s_\ell^*\operatorname{Four}_\psi(\Q_{\mathfrak Z})
 \simeq s_\ell^*o_*\Q_{\mathbb P(\mathcal V^\vee)}[-2t](-t)
 \simeq\iota_*\Q_{\Xi_f}[-2t](-t).
\]
Let
\[
 h:\mathfrak Z\longrightarrow[\mathcal Z/\mathbb G_m],
 \qquad [(z,\xi)]\longmapsto[z].
\]
The equivariant isomorphism
$f^*Rj_{0*}\Q[-1]\simeq
R\operatorname{pr}_{\mathcal Z!}\operatorname{ev}_f^*\mathcal L_\psi$
obtained above descends to
\[
 \bar f^*Rj_*\Q[-1]\simeq Rh_!\ell^*\mathcal L_\psi.
\]
The maps fit into the commutative diagram
\[
 \begin{tikzcd}[column sep=large,row sep=large]
 & \mathfrak Z \arrow[r,"h"] \arrow[d,"q"']
 & { [\mathcal Z/\mathbb G_m] } \arrow[d,"\rho"] \\
 \Xi_f \arrow[r,"\iota"] \arrow[rr,bend right=25,"p"']
 & \mathbb P(\mathcal V^\vee) \arrow[r,"\pi"] & S.
 \end{tikzcd}
\]
Hence
\[
 \begin{aligned}
 R\rho_!\bar f^*Rj_*\Q[-1]
 &\simeq R\rho_!Rh_!\ell^*\mathcal L_\psi
 \simeq R\pi_!Rq_!\ell^*\mathcal L_\psi\\
 &\simeq R\pi_!\iota_*\Q_{\Xi_f}[-2t](-t)
 \simeq Rp_*\Q_{\Xi_f}[-2t](-t),
 \end{aligned}
\]
where $Rp_!=Rp_*$ because $\Xi_f$ is closed in the projective bundle
$\mathbb P(\mathcal V^\vee)$ over $S$.
\end{proof}

\subsection{Proof of the boundary coefficient theorem}\label{bd:sec:lifting}

We now apply Lemma~\ref{bd:lem:fourier} to the lifting problem for $\beta_{U,W}$ in
\eqref{bd:eq:beta}. In the proof below, we will stratify the stack of pairs
\[
 (x,L),\qquad x\in V,\quad L\subset E_x,\quad \dim L=e,
\]
modulo the $\mathbb G_m$-action, according to the dimensions of
$L\cap N$ and $p_x(L)$. We now construct a vector bundle model
for each stratum.
Fix $u\le m$ and $v\le n$ with $u+v=e$, and recall that
\[
 S=S_{u,v}=\Gr_u(M)\times\Gr_v(N)
\]
parameterizes pairs $(U,W)$ of subrepresentations of dimensions $u,v$.
Assume $S$ is nonempty, and write $\mathcal O_S$ for its structure sheaf.
Recall that $V=\Ext^1_Q(M,N)$ and $r=\dim_k V$. Put
\[
 \mathcal V=V\otimes_k\mathcal O_S,\qquad
 t_{u,v}=r+\eu{u}{n-v},
\]
and let $\beta_{U,W}$ be the map in \eqref{bd:eq:beta}.

Let $\mathcal M$ be the representation over $S$ given by
\[
 \mathcal M_i=M_i\otimes_k\mathcal O_S,\qquad
 \mathcal M_h=M_h\otimes\operatorname{id}_{\mathcal O_S},
\]
and define $\mathcal N$ similarly from $N$.
The universal subrepresentations $\mathcal U\subset\mathcal M$ and
$\mathcal W\subset\mathcal N$ satisfy, at $(U,W)\in S$,
\[
 \mathcal U_{(U,W)}=U,\qquad
 \mathcal W_{(U,W)}=W,\qquad
 (\mathcal N/\mathcal W)_{(U,W)}=N/W.
\]
The vector bundles
\[
 \mathcal D^0=\bigoplus_{i\in Q_0}
     \Hom_{\mathcal O_S}(\mathcal U_i,\mathcal N_i/\mathcal W_i),\qquad
 \mathcal D^1=\bigoplus_{h:i\to j,\ h\in Q_1}
     \Hom_{\mathcal O_S}(\mathcal U_i,\mathcal N_j/\mathcal W_j)
\]
collect vertex maps $\lambda_i:U_i\to N_i/W_i$ and arrow maps
$U_i\to N_j/W_j$. Recall the chosen linear section
\[
 \sigma:V\longrightarrow C^1(M,N)
 =\bigoplus_{h:i\to j}\Hom_k(M_i,N_j)
\]
of the quotient onto $V$. Its component $\sigma(x)_h$ is the
upper-right entry of the arrow matrix of $E_x$.
For a fixed pair $(U,W)\in S$, define
\[
 \gamma_{U,W}:V\longrightarrow\mathcal D^1_{(U,W)}
 =\bigoplus_{h:i\to j,\ h\in Q_1}\Hom_k(U_i,N_j/W_j)
\]
by assigning to $x$ the collection of composites
\[
 \gamma_{U,W}(x)_h:\ U_i\hookrightarrow M_i
 \xrightarrow{\sigma(x)_h}N_j\twoheadrightarrow N_j/W_j.
\]
These fibre maps assemble into a bundle map
$\gamma:\mathcal V\to\mathcal D^1$.
Also define $d:\mathcal D^0\to \mathcal D^1$ by
\[
 d(\lambda)_h=(\mathcal N/\mathcal W)_h\lambda_i-\lambda_j\mathcal U_h.
\]

\begin{proposition}[Vector bundle model]\label{bd:prop:lifting-bundle}
There is a vector bundle $\mathcal Z\to S$ of rank $t_{u,v}$
whose fibre at $(U,W)$ parameterizes the pairs
\[
 (x,L),\qquad x\in V,\quad L\subset E_x,\quad
 L\cap N=W,\quad p_x(L)=U.
\]
The map
\[
 f:\mathcal Z\longrightarrow\mathcal V,\qquad
 ((U,W),x,L)\longmapsto((U,W),x),
\]
is a morphism of vector bundles.
On every geometric fibre, there is an exact sequence
\[
 0\longrightarrow\Hom_Q(U,N/W)\longrightarrow\mathcal Z_{(U,W)}
 \xrightarrow{f_{U,W}}\ker\beta_{U,W}\longrightarrow0.
\]
The map $f$ is equivariant for scalar multiplication on both bundles.
\end{proposition}

\begin{proof}
We now describe the pairs $(x,L)$ by linear equations in $(\lambda,x)$.
Fix $(U,W)$ and $x\in V$. Choose subspaces satisfying
\[
 L_i\subset(E_x)_i=N_i\oplus M_i,\qquad
 L_i\cap N_i=W_i,\qquad p_{x,i}(L_i)=U_i.
\]
For each $a\in U_i$, the equality $p_{x,i}(L_i)=U_i$ allows us
to choose $n\in N_i$ with $(n,a)\in L_i$.
If $(n',a)\in L_i$ is another choice, then
$(n-n',0)\in L_i\cap N_i=W_i$, so $n+W_i=n'+W_i$.
Thus $\lambda_i(a)=n+W_i$ is independent of the choice of $n$
and defines a linear map $\lambda_i:U_i\to N_i/W_i$.
By construction,
\[
 L_i/W_i=\{(\lambda_i(a),a):a\in U_i\}
 \subset (N_i/W_i)\oplus M_i.
\]
For an arrow $h:i\to j$, the arrow map of $E_x/W$ sends
\[
 (\lambda_i(a),a)\longmapsto
 \bigl((N/W)_h\lambda_i(a)+\gamma(x)_h(a),\,U_h(a)\bigr),
 \qquad a\in U_i.
\]
This lies in $L_j/W_j$ exactly when its first entry equals
$\lambda_j(U_h(a))$. Thus the subspaces $L_i$ form a
subrepresentation of $E_x$ exactly when,
for every arrow $h:i\to j$,
\[
 (N/W)_h\lambda_i+\gamma(x)_h=\lambda_jU_h.
\]
Equivalently, $d\lambda+\gamma(x)=0$. Define
\[
 \mathcal Z=\ker\bigl((d,\gamma):\mathcal D^0\oplus\mathcal V\to \mathcal D^1\bigr),
 \qquad f:\mathcal Z\longrightarrow\mathcal V,\quad(\lambda,x)\longmapsto x.
\]
This construction also works for families over $S$, identifying
this closed subscheme with the required family of pairs $(x,L)$.

To show that $\mathcal Z$ is a vector bundle, we prove that
$(d,\gamma)$ is surjective on every geometric fibre.
At a geometric point $(U,W)$, write $d_{U,W}$ and $\gamma_{U,W}$
for the fibre maps. Applying the description of extensions in
Subsection~\ref{cc:sec:extension-hall} to $(U,N/W)$, we obtain
\[
 \operatorname{coker}d_{U,W}
 =\mathcal D^1_{(U,W)}/\operatorname{im}d_{U,W}
 \simeq\Ext^1_Q(U,N/W).
\]
The composite
\[
 V\xrightarrow{\gamma_{U,W}}\mathcal D^1_{(U,W)}
 \twoheadrightarrow\operatorname{coker}d_{U,W}
 \xrightarrow{\sim}\Ext^1_Q(U,N/W)
\]
sends $x$ to the quotient extension
$0\to N/W\to p_x^{-1}(U)/W\to U\to0$.
It is therefore $\beta_{U,W}$ from \eqref{bd:eq:beta}.
Since $\beta_{U,W}$ is surjective,
$\operatorname{im}d_{U,W}+\operatorname{im}\gamma_{U,W}=\mathcal D^1_{(U,W)}$.
Thus $(d,\gamma)$ is a surjective bundle map, and its kernel is a
vector bundle in the exact sequence
\begin{equation}\label{bd:eq:Zexact}
 0\longrightarrow\mathcal Z\longrightarrow \mathcal D^0\oplus\mathcal V
 \xrightarrow{(d,\gamma)}\mathcal D^1\longrightarrow0.
\end{equation}
Its rank is
\[
 \operatorname{rk}\mathcal Z
 =r+\operatorname{rk}\mathcal D^0-\operatorname{rk}\mathcal D^1
 =r+\eu{u}{n-v}=t_{u,v}.
\]
Moreover, $f(\lambda,x)=x$ has image $\ker\beta_{U,W}$ on each fibre:
solutions to $d\lambda=-\gamma(x)$ exist exactly when
$\beta_{U,W}(x)=0$, and two solutions differ by an element of
$\ker d_{U,W}=\Hom_Q(U,N/W)$. This gives the asserted exact sequence.

Scaling $(\lambda,x)$ by $(c\lambda,cx)$ makes $f$ equivariant.
Passing to the quotient gives the map to the extension parameter stack:
\[
 \bar f:[\mathcal Z/\mathbb G_m]\longrightarrow[\mathcal V/\mathbb G_m]
 \simeq S\times\mathfrak U_V\longrightarrow\mathfrak U_V.\qedhere
\]
\end{proof}

Recall that the locus $\Xi_f$ in Lemma~\ref{bd:lem:fourier} is defined by
$f_{U,W}^\vee(\eta)=0$. To identify this locus, we consider the dual bundle map
$f^\vee:\mathcal V^\vee\to\mathcal Z^\vee$ of the map constructed in
Proposition~\ref{bd:prop:lifting-bundle}.
Write $f_{U,W}^\vee:D_kV\to\mathcal Z_{(U,W)}^\vee$ for its fibre map
at $(U,W)\in S$, and let
$\pi:S\times\mathbb P(D_kV)=\mathbb P(\mathcal V^\vee)\to S$ be the projection.

\begin{lemma}[Identification of the dual zero locus]\label{bd:lem:dual-zero}
On every geometric fibre,
\[
 \ker f_{U,W}^\vee
 =\operatorname{im}(D_k\beta_{U,W})
 =\{\eta\in D_kV:\eta|_{\ker\beta_{U,W}}=0\}.
\]
The locus
\[
 \Xi_{u,v}=\{(U,W,[\eta])\in S_{u,v}\times\mathbb P(D_kV):
                 \eta|_{\ker\beta_{U,W}}=0\}
\]
has a closed subscheme structure given by
\[
 \Xi_{u,v}:=\operatorname{Zero}\!\left(
 \mathcal O(-1)\hookrightarrow\pi^*\mathcal V^\vee
 \xrightarrow{\pi^*f^\vee}\pi^*\mathcal Z^\vee\right).
\]
With this structure, $\Xi_{u,v}=\Xi_f$.
\end{lemma}

\begin{proof}
By Proposition~\ref{bd:prop:lifting-bundle}, there is an exact sequence
\[
 0\longrightarrow\Hom_Q(U,N/W)\longrightarrow\mathcal Z_{(U,W)}
 \xrightarrow{f_{U,W}}\ker\beta_{U,W}\longrightarrow0.
\]
For $\eta\in D_kV$ and $(\lambda,x)\in\mathcal Z_{(U,W)}$, we have
\[
 \bigl(f_{U,W}^\vee(\eta)\bigr)(\lambda,x)
 =\eta\bigl(f_{U,W}(\lambda,x)\bigr)=\eta(x).
\]
Since $\operatorname{im}f_{U,W}=\ker\beta_{U,W}$, this functional vanishes on
$\mathcal Z_{(U,W)}$ if and only if $\eta$ vanishes on $\ker\beta_{U,W}$.
Hence
\[
 \ker f_{U,W}^\vee
 =\{\eta\in D_kV:\eta|_{\ker\beta_{U,W}}=0\}
 =\operatorname{im}(D_k\beta_{U,W}).
\]
The last equality holds because $\beta_{U,W}$ is surjective:
a functional on $V$ vanishes on $\ker\beta_{U,W}$ exactly when
it factors through $\beta_{U,W}$.

At $(U,W,[\eta])$, the tautological line $\mathcal O(-1)$ is spanned by
$\eta$, so the displayed bundle map vanishes precisely when
$f_{U,W}^\vee(\eta)=0$, or equivalently, when $\eta|_{\ker\beta_{U,W}}=0$.
This is the condition in \eqref{bd:eq:xi}. Since the same bundle map defines
$\Xi_f$, we obtain $\Xi_{u,v}=\Xi_f$ as closed subschemes.
\end{proof}

\begin{proof}[Proof of Theorem~\ref{bd:thm:main}]
Form the Cartesian square
\[
 \begin{tikzcd}
 \mathfrak H_e=\mathfrak U_V\times_{\M_{m+n}}\mathfrak G_{m+n,e}
   \arrow[r]\arrow[d,"\widetilde\rho_{\alpha,e}"']
   &\mathfrak G_{m+n,e}\arrow[d,"\rho_{\alpha,e}"]\\
 \mathfrak U_V\arrow[r,"\phi"']&\M_{m+n}.
 \end{tikzcd}
\]
Recall that $\mathcal B_{M,N}=\phi_!Rj_*\Q[-1]$ and
$\mathcal K_{(\alpha-e,e)}\simeq R\rho_{\alpha,e!}\Q$.
Proper base change in the Cartesian square gives
\[
 \phi^*\mathcal K_{(\alpha-e,e)}
 \simeq\phi^*R\rho_{\alpha,e!}\Q
 \simeq R\widetilde\rho_{\alpha,e!}\Q_{\mathfrak H_e}.
\]
Using the definition of the pairing and the projection formula for
$\phi$ and $\widetilde\rho_{\alpha,e}$, we obtain
\[
 \begin{aligned}
 \mathsf P(\mathcal B_{M,N},\mathcal K_{(\alpha-e,e)})
 &=R\Gamma_c(\M_\alpha,(\phi_!Rj_*\Q[-1])\otimes\mathcal K_{(\alpha-e,e)})\\
 &\simeq R\Gamma_c(\mathfrak U_V,Rj_*\Q[-1]\otimes\phi^*\mathcal K_{(\alpha-e,e)})\\
 &\simeq R\Gamma_c(\mathfrak U_V,Rj_*\Q[-1]\otimes R\widetilde\rho_{\alpha,e!}\Q_{\mathfrak H_e})\\
 &\simeq R\Gamma_c(\mathfrak U_V,R\widetilde\rho_{\alpha,e!}
       (\widetilde\rho_{\alpha,e}^*Rj_*\Q[-1]))\\
 &\simeq R\Gamma_c(\mathfrak H_e,\widetilde\rho_{\alpha,e}^*Rj_*\Q[-1]).
 \end{aligned}
\]
Stratify $\mathfrak H_e$ by the locally closed loci $\mathfrak H_{u,v}$
defined by
\[
 \dim p_x(L)=u,\qquad \dim(L\cap N)=v,\qquad u+v=e.
\]
Proposition~\ref{bd:prop:lifting-bundle} identifies $\mathfrak H_{u,v}$ with
$[\mathcal Z/\mathbb G_m]$ over $S_{u,v}$.
The conditions $\dim(L\cap N)\ge v$ are closed, so an ordering of
the finitely many strata compatible with this partial order gives a
localization filtration. Its layers are
\[
 R\Gamma_c([\mathcal Z/\mathbb G_m],\bar f^*Rj_*\Q[-1])
 \simeq R\Gamma_c(\Xi_{u,v},\Q)[-2t_{u,v}](-t_{u,v}),
\]
where the isomorphism follows from Lemma~\ref{bd:lem:fourier} over $S_{u,v}$
and the identification $\Xi_f=\Xi_{u,v}$ in Lemma~\ref{bd:lem:dual-zero}.
This proves \eqref{bd:eq:layers}. Moreover, $t_{u,v}=\operatorname{rk}\mathcal Z\ge0$.
Taking traces gives
\[
 b_e(q)=\sum_{u+v=e}q^{t_{u,v}}
       \Tr(\Fr,R\Gamma_c(\Xi_{u,v},\Q))
       =\sum_{u+v=e}q^{t_{u,v}}\#\Xi_{u,v}(k),
\]
by additivity and the trace formula. All constructions commute with
finite extension of $k$.
\end{proof}

\subsection{The boundary character}\label{gm:sec:dual-data}

We describe the fibres of
$\Xi_{u,v}\to\mathbb P(D_k\Ext^1_Q(M,N))$
and compute their contributions to the boundary character.
Recall that $\Lambda$ is the skew-symmetric form defined in \eqref{cc:skew-form}
and that $X^zX^{z'}=w^{\Lambda(z,z')}X^{z+z'}$.
Applying the normalized character construction
\eqref{geo:character} to $\mathcal B_{M,N}$, we obtain
\begin{equation}\label{gm:boundary-polynomial}
 {\widehat{\mathcal B}_{M,N}
 =\Phi([\mathcal B_{M,N}])
 =\sum_e w^{-\eu{e}{m+n-e}}b_e(q)X^{(-E(m+n)-Be,m+n-e)},}
\end{equation}
where $m=\dim M$, $n=\dim N$, and $\Phi$ is as in
\eqref{cc:Hall-character-normalization}.
Let $\Xi_{\eta;u,v}$ be the fibre of $\Xi_{u,v}\to\mathbb P(D_kV)$
at $[\eta]$. Theorem~\ref{bd:thm:main} gives
\[
 \widehat{\mathcal B}_{M,N}
 =\sum_{[\eta]\in\mathbb P(D_kV)(k)}\widehat{\mathcal B}_\eta,
\]
with
\begin{equation}\label{gm:boundary-fibre-sum}
 \begin{aligned}
 \widehat{\mathcal B}_\eta
 =\sum_{u,v}{}&w^{2(r+\eu{u}{n-v})}w^{-\eu{u+v}{m+n-u-v}}\\
 &\cdot|\Xi_{\eta;u,v}(k)|X^{(-E(m+n)-B(u+v),m+n-u-v)}.
 \end{aligned}
\end{equation}
The shift $[-2t_{u,v}](-t_{u,v})$ in \eqref{bd:eq:layers}
contributes the trace factor $q^{t_{u,v}}=w^{2t_{u,v}}$,
where $t_{u,v}=r+\eu{u}{n-v}$. This is the first power of $w$.
The second is the character normalization in total degree $u+v$.

The category $\operatorname{rep}_kQ$ is hereditary.
The AR formulas give natural isomorphisms
\begin{equation}\label{gm:AR-adjunction}
 D_k\Ext^1_Q(M,N)
 \simeq\Hom_Q(N,\tau M)
 \simeq\Hom_Q(\tau^-N,M).
\end{equation}
Here $\tau$ denotes the AR translate on $\operatorname{rep}_kQ$,
and $\tau^-$ its inverse translate.

For $0\ne\eta\in D_kV$, denote the corresponding maps by
$\eta:N\to\tau M$ and $\eta^\flat:\tau^-N\to M$. Set
\begin{equation}\label{gm:reverse-data}
 K_\eta=\operatorname{im}\eta^\flat,\qquad
 J_\eta=\operatorname{coker}\eta^\flat,\qquad
 D_\eta=\ker\eta.
\end{equation}
Write
\[
 \kappa=\dim K_\eta,\qquad j=\dim J_\eta,\qquad d=\dim D_\eta.
\]

\begin{lemma}[Fibres of the dual zero locus]\label{gm:dual-fibre}
For every $[\eta]\in\mathbb P(D_kV)$,
\begin{equation}\label{gm:dual-fibre-formula}
 (\Xi_{\eta;u,v})_{\rm red}
 \simeq
 \bigl(\Gr_{u-\kappa}(J_\eta)\times\Gr_v(D_\eta)\bigr)_{\rm red}.
\end{equation}
\end{lemma}

\begin{proof}
By Lemma~\ref{bd:lem:dual-zero}, a geometric point
$(U,W)\in S_{u,v}$ lies in $\Xi_{\eta;u,v}$ if and only if
$\eta\in\operatorname{im}D_k\beta_{U,W}$.
Recall $\beta_{U,W}=b_W\circ a_U$ from \eqref{bd:eq:beta}.
Under the AR isomorphisms \eqref{gm:AR-adjunction}, its dual
$D_k\beta_{U,W}=D_ka_U\circ D_kb_W$ is the composite
\[
 \Hom_Q(N/W,\tau U)\xrightarrow{D_kb_W}
 \Hom_Q(N,\tau U)\xrightarrow{D_ka_U}\Hom_Q(N,\tau M).
\]
Explicitly, for $h:N/W\to\tau U$ and $g:N\to\tau U$,
\[
 D_kb_W(h)=\bigl(N\twoheadrightarrow N/W\xrightarrow{h}\tau U\bigr),
 \qquad
 D_ka_U(g)=\bigl(N\xrightarrow{g}\tau U\longrightarrow\tau M\bigr).
\]
The map $\tau U\to\tau M$, induced by $U\hookrightarrow M$, is injective.
Indeed, applying $\Hom_Q(-,kQ)$ to
$0\to U\to M\to M/U\to0$ gives
\[
 \Ext^1_Q(M,kQ)\longrightarrow\Ext^1_Q(U,kQ)
 \longrightarrow\Ext^2_Q(M/U,kQ)=0,
\]
where the last equality follows from heredity.
The natural isomorphism $\tau X\simeq D_k\Ext^1_Q(X,kQ)$ identifies
the dual of this surjection with
\[
 \tau U\simeq D_k\Ext^1_Q(U,kQ)
 \hookrightarrow D_k\Ext^1_Q(M,kQ)\simeq\tau M.
\]
By naturality, this is the map induced by $U\hookrightarrow M$.
If $\eta$ factors through this injection via
$\widetilde\eta:N\to\tau U$, injectivity shows that
$\widetilde\eta|_W=0$ if and only if $\eta|_W=0$.
Hence $\widetilde\eta$ descends to a map $h:N/W\to\tau U$
exactly when $W\subseteq\ker\eta=D_\eta$.

By naturality of the other AR identification in \eqref{gm:AR-adjunction},
$\eta$ factors through $\tau U\hookrightarrow\tau M$ if and only if
$\eta^\flat$ factors through $U\hookrightarrow M$.
The latter condition is exactly
$\operatorname{im}\eta^\flat=K_\eta\subseteq U$.
Consequently, the condition $\eta\in\operatorname{im}D_k\beta_{U,W}$ is
equivalent to the two inclusions
\[
 W\subseteq D_\eta,\qquad K_\eta\subseteq U.
\]
A subrepresentation
$U\subseteq M$ containing $K_\eta$ determines
a subrepresentation $U/K_\eta\subseteq\operatorname{coker}\eta^\flat=J_\eta$
of dimension $u-\kappa$.
Conversely, taking the inverse image under $M\twoheadrightarrow J_\eta$
recovers $U$. The map
\[
 (U,W)\longmapsto(U/K_\eta,W)
\]
and this inverse-image construction identify the closed subscheme of
$S_{u,v}$ defined by $K_\eta\subseteq U$ and $W\subseteq D_\eta$
with $\Gr_{u-\kappa}(J_\eta)\times\Gr_v(D_\eta)$.
This closed subscheme and $\Xi_{\eta;u,v}$ have the same geometric
points, so their reductions in $S_{u,v}$ coincide. This proves
\eqref{gm:dual-fibre-formula}.
\end{proof}

\begin{remark}\label{gm:reduction}
For a finite-type $k$-scheme $X$, pullback along $X_{\rm red}\to X$
induces an equivalence between the categories of \'etale sheaves on
$X$ and $X_{\rm red}$ by \cite[Tag~03SI]{StacksProject}.
In particular, constant $\ell$-adic
complexes and their compactly supported cohomology are unchanged by
reduction. Moreover, $X_{\rm red}(k')=X(k')$ for every finite extension
$k'/k$, since $\Spec k'$ is reduced. Thus the reductions in
Lemma~\ref{gm:dual-fibre} do not affect the cohomology or point counts
used below.
\end{remark}

Taking point counts in Lemma~\ref{gm:dual-fibre} and setting
$u=\kappa+x$, $v=y$ gives the following expression directly
from the fibres:
\begin{equation}\label{gm:grassmannian-boundary}
 \begin{aligned}
 \widehat{\mathcal B}_\eta
 =\sum_{x,y}{}&w^{2(r+\eu{\kappa+x}{n-y})}w^{-\eu{\kappa+x+y}{m+n-\kappa-x-y}}\\
 &\cdot|\Gr_x(J_\eta)(k)|\,|\Gr_y(D_\eta)(k)|
 X^{(-E(m+n)-B(\kappa+x+y),\,m+n-\kappa-x-y)}.
 \end{aligned}
\end{equation}
Here $x=\dim(U/K_\eta)$ and $y=\dim W$ are the dimensions of
the corresponding subrepresentations of $J_\eta$ and $D_\eta$.
We now compare these exponents with those in $X_{J_\eta}$ and
$X_{D_\eta}$. Since $m=\kappa+j$ and $B=E^T-E$,
the boundary exponent minus the sum of the two character exponents
is $(p_\eta,n-d)$, where
\begin{equation}\label{gm:projective-vector}
 p_\eta:=-E(n-d)-E^T\kappa.
\end{equation}
Equivalently,
\begin{equation}\label{gm:monomial-match}
 -E(m+n)-B\kappa=-E(j+d)+p_\eta.
\end{equation}
For the frozen component, the boundary term has exponent
$m+n-\kappa-x-y$, whereas the corresponding product of terms in
$X_{J_\eta}$ and $X_{D_\eta}$ has exponent $j+d-x-y$.
Their difference is therefore $n-d$.

\begin{proposition}\label{gm:ordered}
For every $[\eta]\in\mathbb P(D_kV)(k)$,
\begin{equation}\label{gm:ordered-formula}
 \widehat{\mathcal B}_\eta
 =w^{2r+\eu{\kappa}{n-d}}
 X_{J_\eta}X^{(p_\eta,n-d)}X_{D_\eta}.
\end{equation}
\end{proposition}

\begin{proof}
Put $z=(p_\eta,n-d)$. The weight exponent in
\eqref{gm:grassmannian-boundary} is
\[
 H(x,y)=2r+2\eu{\kappa+x}{n-y}
             -\eu{\kappa+x+y}{m+n-\kappa-x-y}.
\]
Since $m=\kappa+j$, expansion using
\eqref{gm:projective-vector} gives
\begin{align*}
 H(x,y)={}&2r+\eu{\kappa}{n-j}
              -\eu{x}{j-x}-\eu{y}{d-y}\\
 &+\eu{x}{d}-\eu{y}{j}+\eu{y}{x}-\eu{x}{y}
              -p_\eta^Tx+p_\eta^Ty.
\end{align*}
On the other hand, the exponent contributed by the three-factor
product is
\[
 \begin{aligned}
 &\Lambda((-Ej-Bx,j-x),z)
   +\Lambda((-Ej-Bx,j-x),(-Ed-By,d-y))\\
 &\qquad+\Lambda(z,(-Ed-By,d-y))\\
 &\quad=\eu{\kappa}{d-j}
      +\eu{x}{d}-\eu{y}{j}+\eu{y}{x}-\eu{x}{y}
      -p_\eta^Tx+p_\eta^Ty.
 \end{aligned}
\]
The character weights supply
$-\eu{x}{j-x}-\eu{y}{d-y}$. The remaining exponent is
$2r+\eu{\kappa}{n-j}-\eu{\kappa}{d-j}
=2r+\eu{\kappa}{n-d}$.
Together with \eqref{gm:monomial-match}, this proves the formula.
\end{proof}

\subsection{The general multiplication formula}

The extension triangle first gives a formula with its boundary
character. Substituting Proposition~\ref{gm:ordered} then gives the
explicit expression below. Compare
\cite[Theorem~7.4 and Corollary~7.5]{ChenDingZhang2023}.

\begin{theorem}[Multiplication from the extension triangle]
\label{gm:general-multiplication}
Let $M,N$ be representations with $r=\dim\Ext^1_Q(M,N)>0$. Then
\begin{equation}\label{gm:boundary-product}
 [r]_qX_MX_N
 =\sum_{[x]\in\mathbb P(V)(k)}X_{E_x}
         +w^{-2r}\widehat{\mathcal B}_{M,N}.
\end{equation}
Explicitly,
\begin{equation}\label{gm:general-product}
 {[r]_qX_MX_N
 =\sum_{[x]\in\mathbb P(V)(k)}X_{E_x}
 +\sum_{[\eta]\in\mathbb P(D_kV)(k)}
 w^{\eu{\kappa}{n-d}}X_{J_\eta}X^{(p_\eta,n-d)}X_{D_\eta}.}
\end{equation}
\end{theorem}

\begin{proof}
Apply $\Phi$ to the triangle and the filtration in
Theorem~\ref{cc:sheaf-exchange}. The trace formula and
\eqref{cc:Hall-character-normalization} give
\[
 \Phi([\mathcal E_{M,N}])
 =\sum_{[x]\in\mathbb P(V)(k)}X_{E_x},\qquad
 \Phi([A_M\star A_N])=X_MX_N.
\]
The middle filtration contributes $[r]_q$, and the shift $[2r](r)$
of the boundary term contributes $w^{-2r}$. This proves
\eqref{gm:boundary-product}. Substituting \eqref{gm:ordered-formula}
gives \eqref{gm:general-product}.
\end{proof}

\subsection{Generic exchange relations}
\label{gm:sec:exchange}

We specialize Theorem~\ref{gm:general-multiplication} to $r=1$.
Choose nonzero $x\in V$ and $\eta\in D_kV$,
and write $E_x$ for the middle term of
\[
 0\longrightarrow N\longrightarrow E_x\longrightarrow M\longrightarrow0.
\]
Since $[1]_q=1$ and both projective spaces consist of one point,
\eqref{gm:boundary-product} has two terms on the right.
We use rigidity and counting polynomials to pass to the formal parameter
$w$ and compare this formula with quantum exchange relations below.

We first record the representation-theoretic meaning of $p_\eta$.
The following description
holds for arbitrary $r$. Write
$\nu=D_k\Hom_Q(-,kQ)$ for the Nakayama equivalence, and decompose
\[
 M=M_0\oplus P,\qquad\operatorname{coker}\eta=R_\eta\oplus I,
\]
where $P$ is projective, $I$ is injective, and $M_0,R_\eta$ have no
projective and injective summands, respectively.

\begin{lemma}[Projective multiplicities]\label{gm:projective-multiplicities}
The vector $p_\eta$ in \eqref{gm:projective-vector} satisfies
\[
 p_\eta=E\dim I\in\mathbb N^{Q_0}
\]
and is the multiplicity vector of the indecomposable projective
summands of $\nu^{-1}I$.
\end{lemma}

\begin{proof}
Since $P$ is projective, $\Ext^1_Q(P,N)=0$.
Applying \eqref{gm:AR-adjunction} to $(P,N)$ gives
$\Hom_Q(\tau^-N,P)\simeq D_k\Ext^1_Q(P,N)=0$,
so the $P$-component of $\eta^\flat:\tau^-N\to M_0\oplus P$
vanishes. Hence $\operatorname{im}\eta^\flat\subseteq M_0$.
Applying the right exact functor $\tau^-$ to
$N\xrightarrow{\eta}\tau M_0\to R_\eta\oplus I\to0$ gives
\[
 \tau^-N\xrightarrow{\eta^\flat}M_0\longrightarrow\tau^-R_\eta\longrightarrow0,
\]
using $\tau^-\tau M_0\simeq M_0$ and $\tau^-I=0$.
Since $K_\eta=\operatorname{im}\eta^\flat$, this gives the short exact sequence
\[
 0\longrightarrow K_\eta\longrightarrow M_0
 \longrightarrow\tau^-R_\eta\longrightarrow0.
\]
Thus
\[
 J_\eta\simeq\tau^-R_\eta\oplus P,\qquad
 \kappa=\dim M_0-\dim\tau^-R_\eta.
\]
The exact sequence
\[
 0\longrightarrow D_\eta\longrightarrow N\xrightarrow{\eta}\tau M_0
 \longrightarrow R_\eta\oplus I\longrightarrow0
\]
gives $\dim N-d=\dim\tau M_0-\dim R_\eta-\dim I$.
The Coxeter identities
\[
 E\dim\tau M_0=-E^T\dim M_0,\qquad
 E\dim R_\eta=-E^T\dim\tau^-R_\eta
\]
therefore yield
\[
 \begin{aligned}
 E(\dim N-d)
 &=-E^T\dim M_0+E^T\dim\tau^-R_\eta-E\dim I\\
 &=-E^T\kappa-E\dim I.
 \end{aligned}
\]
Thus $p_\eta=E\dim I$.
Since $I$ is injective, for each vertex $i$ we have
\[
 (p_\eta)_i=\eu{\dim S_i}{\dim I}
 =\dim_k\Hom_Q(S_i,I)\in\mathbb N.
\]
By the Nakayama equivalence, this is the multiplicity of the
projective cover of $S_i$ in $\nu^{-1}I$.
\end{proof}

We use the generic characters of Subsection~\ref{sec:cluster}, expressed
in the fixed initial quantum torus with formal parameter $w$.

\begin{corollary}[Quantum exchange relation]\label{gm:exchange-mutation}
Let $T=\overline T\oplus M$ and
$T'=\overline T\oplus N=\mu_k(T)$ be adjacent reachable support tilting
representations of the same support, with
$\dim_k\Ext^1_Q(M,N)=1$.
{The multiplication formula of Theorem~\ref{gm:general-multiplication}
gives the following quantum exchange relation for the generic characters
in direction $k$ over $\mathbb Z[w^{\pm1}]$.}
\begin{equation}\label{gm:rigid-exchange}
 X_MX_N=X_{E_x}
   +w^{-1-\eu{\kappa}{j+d}}
       X^{(p_\eta,n-d)}X_{J_\eta\oplus D_\eta}.
\end{equation}
\end{corollary}

{For $r=1$, the exchange sequences of \cite[Propositions~2.9--2.10]{Rupel2015}
and rigidity reduce \eqref{gm:general-product} to \eqref{gm:rigid-exchange}.
We use $\eu{\kappa}{n}=-1$ in the exponent.
Here $\eu{\kappa}{j}=0$ and $\Hom_Q(J_\eta,I)=0$ give the required commutation
by Lemma~\ref{gm:projective-multiplicities}. Passing to counting polynomials gives
the generic exchange relation at $k$ by \cite[Theorem~6.3]{Rupel2015}.}

\begin{remark}
The correction monomial in \eqref{gm:rigid-exchange} is visible in the grading
$\deg X^{(a,b)}=a-Bb$, $\deg w=0$.
Indeed, $\deg X_M=-E^Tm$ and
$\deg X_{J_\eta\oplus D_\eta}=-E^T(j+d)$, and
\eqref{gm:monomial-match} gives
\[
 \deg X_M+\deg X_N=\deg X_{J_\eta\oplus D_\eta}+p_\eta-B(n-d).
\]
Thus $X^{(p_\eta,n-d)}$ supplies the missing degree.
\end{remark}

\section*{AI Declaration}
We had developed the first part of the paper before using generative AI on this project. We used AI in studying the boundary term. More precisely, Chao Shen consulted OpenAI's GPT-6 Astra, which suggested a standard Fourier transform
argument and helped translate the relevant passages of Laumon's paper from French.
Shen then worked out the calculation in the present setting, checked the details,
and wrote the proof. AI is also used to check references and revise the exposition.
The authors take full responsibility for the manuscript.

\end{document}